\documentclass[12pt,reqno]{amsart}
\usepackage[margin=1in]{geometry}
\usepackage{amsmath,amssymb,amsthm,mathtools,mathrsfs}
\usepackage{xcolor}
\usepackage[hyperfootnotes=false, colorlinks, linkcolor={blue}, citecolor={magenta}, filecolor={blue}, urlcolor={blue}]{hyperref}
\usepackage{etoolbox}
\usepackage{orcidlink}

\numberwithin{equation}{section}

\allowdisplaybreaks
\newcommand{\A}{\mathbb A}
\newcommand{\C}{\mathbb C}
\newcommand{\Q}{\mathbb Q}
\newcommand{\R}{\mathbb R}
\newcommand{\Z}{\mathbb Z}

\newcommand{\Ad}{\operatorname{Ad}}
\newcommand{\diag}{\operatorname{diag}}

\newcommand{\Sp}{\mathrm{Sp}}
\newcommand{\GSp}{\mathrm{GSp}}
\newcommand{\PGSp}{\mathrm{PGSp}}
\newcommand{\GL}{\mathrm{GL}}
\newcommand{\PGL}{\mathrm{PGL}}
\newcommand{\SO}{\mathrm{SO}}

\newcommand{\SL}{\mathrm{SL}}

\newcommand{\Sym}{\operatorname{Sym}}
\newcommand{\tr}{\operatorname{tr}}
\newcommand{\vol}{\operatorname{vol}}

\newcommand{\cD}{\mathcal D}

\newcommand{\Tscr}{\mathcal T}

\newcommand{\p}{\mathfrak p}
\newcommand{\OF}{{\mathfrak o}}
\newcommand{\triv}{1}
\newcommand{\AI}{{\mathcal{AI}}}

\newtheorem{theorem}{Theorem}[section]
\newtheorem{proposition}[theorem]{Proposition}
\newtheorem{lemma}[theorem]{Lemma}
\newtheorem{corollary}[theorem]{Corollary}
\theoremstyle{definition}
\newtheorem{definition}[theorem]{Definition}
\theoremstyle{remark}
\newtheorem{remark}[theorem]{Remark}

\DeclareMathOperator{\Cl}{Cl}

\makeatletter
\@namedef{subjclassname@2020}{
  \textup{2020} Mathematics Subject Classification}
\makeatother

\makeatletter
\patchcmd{\@settitle}
  {\uppercasenonmath\@title}
  {}{}{}
\makeatother

\begin{document}

\title{When GL(2) cusp forms meet CM points in GSp(4)}
\author{Shenghao Hua~\orcidlink{0000-0002-7210-2650}}
\address{Shanghai Institute for Mathematics and Interdisciplinary Sciences (SIMIS), Shanghai 200433, China}
\address{Research Institute of Intelligent Complex Systems, Fudan University, Shanghai, 200433, China}
\email{huashenghao@vip.qq.com}

\author{Shaoyun Yi~\orcidlink{0000-0002-1251-8611}}
\address{School of Mathematical Sciences, Xiamen University, Xiamen, Fujian 361005, China}
\email{yishaoyun926@xmu.edu.cn}

\date{}

\begin{abstract}
Duke's equidistribution theorem shows that the nonconstant contribution of the $\mathrm{GL}(2)$ automorphic spectrum to CM packets becomes asymptotically negligible.
CM packets on the Levi base index the fibers of Bessel quotients in $\mathrm{GSp}(4)$, and the symmetric cube correspondence carries $\mathrm{GL}(2)$ cusp forms to $\mathrm{GL}(4)$.
In the holomorphic setting considered here, these two constructions meet inside $\mathrm{GSp}(4)$.
We prove averaged decay for the resulting class group periods over imaginary quadratic fields, expressed as averages of conditional fiber integrals.
More precisely, let $f$ be a non-CM holomorphic Hecke newform with trivial central character.
The mean absolute value of the associated class-group averages tends to zero along the transverse unipotent direction.

\end{abstract}


\maketitle

\section{Introduction and statements}
The distribution of automorphic forms over class groups has been an interesting topic in analytic number theory.
For a negative fundamental discriminant $d$, the ideal classes of $K_d=\Q(\sqrt d)$ give rise to Heegner points on modular curves, and therefore provide arithmetic points at which automorphic forms can be evaluated. A fundamental question is to understand the limiting behavior of these values as the discriminant $|d|$ tends to infinity.
Duke's theorem~\cite{Duke1988} establishes that the corresponding Heegner points become equidistributed with respect to the hyperbolic measure.
In particular, the class-group average of a fixed weight-zero cuspidal scalar observable on the modular curve tends to zero.
Duke's theorem reveals a deep connection between class groups, homogeneous dynamics, and automorphic forms, and motivates further studies of toric periods and automorphic $L$-functions.

The adelic interpretation of CM packets forms part of the broader study of arithmetic torus orbits on homogeneous spaces.
Michel and Venkatesh clarify the interaction among equidistribution, automorphic $L$-functions,
and ergodic methods in this setting \cite{MichelVenkatesh2006}, while
Einsiedler, Lindenstrauss, Michel, and Venkatesh analyze the distribution of periodic torus orbits on higher-dimensional homogeneous spaces
\cite{EinsiedlerLindenstraussMichelVenkatesh2009}.
These results provide the context for our adelic formulation, in which the observable carries an additional nontrivial Fourier character along a unipotent fiber.

In the holomorphic setting, CM tori enter Bessel periods on $\GSp(4)$ for a representation whose spin transfer is the symmetric cube of a $\GL(2)$ cusp form.

Fix a non-CM holomorphic Hecke newform $f$ of even weight $\kappa\geq2$
and trivial central character, with associated representation $\pi_f$.
Let $\Pi$ be a cuspidal realization of its symmetric-cube spin transfer,
with $\Pi_\infty\simeq B_{\kappa+1,\kappa-2}$ as in
\cite{Schmidt2017}.
Fix a nonzero factorizable vector $\varphi\in\Pi$ whose real component lies in the holomorphic lowest $K$-type of highest weight $(2\kappa-1,\kappa+1)$. The precise realization, temperedness, and refined special Bessel period formula assumptions are recorded in Subsection~\ref{subsec:standing-data} and remain in force throughout.

For $d<0$ fundamental, let $\chi_d$ be the quadratic character of $K_d$ and $R_d=T_dN$ the Bessel subgroup attached to $S_d$ in
\eqref{eq:Sd}.
Section~\ref{sec:geometry} constructs its compact quotient $[R_d]_Z$, the projection to the CM packet on the Levi base, and the scalar map $\iota_d$ induced by $X\mapsto\tr(S_dX)$; see \eqref{e_iota_d}.
With Tamagawa measure $dr$ of total volume $2$, define our CM-average measure $\mathscr A_d(\varphi)=(\iota_d)_*(\varphi(r)\,dr/2)$.
Its Fourier coefficient is $\mathcal C_{d,b}(\varphi)=\mathcal B_{d,b}(\varphi)/2$,
where $\mathcal B_{d,b}$ is defined in \eqref{eq:Bessel-period}. For our purpose, we estimate the associated Bessel Fourier coefficients as $|d|\to\infty$ and their dyadic means as $D\to\infty$.

Let $M_\Pi$ be the conductor of $\Sym^3\pi_f$ and $M_0=\mathrm{lcm}(8,M_\Pi)$.
Choose a finite set $\Sigma_\varphi$ containing $2$ and all primes where $\Pi_p$ is
ramified or $\varphi_p$ is nonspherical. Let $M_\varphi$ be divisible
by $M_0$ and sufficiently large powers of primes in $\Sigma_\varphi$
so that the bad local data are constant up to integral conjugacy;
see Section~\ref{sec:local}. Fix $a_0\equiv1\pmod4$ coprime to
$M_\varphi$, and put
\[
 \mathcal D(D,a_0)=
 \left\{d<0\ \middle|\,
 \begin{array}{l}
 d\text{ is a fundamental discriminant},\\
 D\leq |d|\leq2D,\quad d\equiv a_0\pmod{M_\varphi}
 \end{array}
 \right\}.
\]

Let $\mathscr P_d^\natural(\varphi)$ be the distribution on
$\Q\backslash\A$ with normalized Fourier coefficients
$\mathcal B^\natural_{d,b}(\varphi)$ defined in
\eqref{eq:normalized-period}. We now state our main result, which gives a quantitative decay estimate for $\mathscr P_d^\natural(\varphi)$.

\begin{theorem}
\label{thm:pointwise-distribution}
Assume the Generalized Riemann Hypothesis (GRH) for
$L(s,\chi_d)$ and $L(s,\Sym^3\pi_f\otimes\chi_d)$
for the negative fundamental discriminants $d\equiv a_0\pmod{M_\varphi}$.
For every $\varepsilon>0$ and every finite-level test function $\Phi$ on $\Q\backslash\A$ that is smooth at the real place, we have
\[
 |\langle\mathscr P_d^\natural(\varphi),\Phi\rangle|
 \ll_{\varphi,\Phi,\varepsilon}|d|^{-1/4+\varepsilon}.
\]
In particular, the pairing tends to zero as $|d|\to\infty$ along this congruence family.
\end{theorem}

Theorem~\ref{thm:distribution} strengthens this to a third-order
seminorm bound uniform over finite test levels. Under the additional mean hypotheses stated there, it also gives
\[
 \frac1D\sum_{d\in\mathcal D(D,a_0)}
 |\langle\mathscr P_d^\natural(\varphi),\Phi\rangle|
 \ll_{\varphi,\Phi,\varepsilon}
 D^{-1/4}(\log D)^{-1/8+\varepsilon}.
\]
For the original measure, Theorem~\ref{thm:raw-exponential} gives
unconditionally
\[
|\langle\mathscr A_d(\varphi),\Phi\rangle|
 \ll_{\varphi,\varepsilon}
 \|\Phi\|_\infty |d|^{\kappa-1/2+\varepsilon}
 \exp\!\left(-\frac{\pi|d|}{2Q_\varphi}\right),
\]
where $Q_\varphi$ is the denominator bound of Lemma~\ref{lem:denominator}.
This estimate applies to every bounded Borel test function $\Phi$ through its supremum norm. The normalized distribution removes this real exponential factor and retains the arithmetic information in the central $L$-values.

Hua--Miao~\cite[Theorems~1.1 and~3.2]{HuaMiao2025} obtain
mean bounds for symmetric-cube Bessel periods from quadratic-twist
moment estimates under GRH.
We combine their moment estimates with local bounds uniform in
$b,d$ to obtain the dyadic mean estimate for
$\mathscr P_d^\natural(\varphi)$.

For $b=m/n$ in lowest terms, $n>0$, define its height by
\begin{equation}\label{e_height_b}
 H(b)=\max\{|m|,n\},\quad H(0)=1. 
\end{equation}
Let $\langle\varphi,\varphi\rangle$ be the squared Tamagawa
Petersson norm, and let $\mathfrak c_d=w(K_d)/h(K_d)$, where $h(K_d)$
is the class number and $w(K_d)$ counts roots of unity. Set
\begin{equation}\label{eq:L-weight}
 \mathscr L_d(\Pi)=
 \frac{\mathfrak c_d}{S_\Pi}
 \frac{\zeta_\Q(2)\zeta_\Q(4)
 |L(1/2,\Pi)L(1/2,\Pi\otimes\chi_d)|}
 {|L(1,\Pi,\Ad)L(1,\chi_d)|}.
\end{equation}
Here $S_\Pi=2^{\ell(\Pi)}$ is the component-group constant of
\cite[Theorem~1]{FurusawaMorimoto2021}, where $\ell(\Pi)$ counts the
cuspidal isobaric summands of the weak spin transfer; here $\ell(\Pi)=1$.
We write $L(s,\Pi)=L(s,\Sym^3\pi_f)$ for the spin $L$-function,
with the same convention for twists.
Automorphic $L$-functions denote their finite parts, i.e., Euler products over all finite primes, whereas
$\zeta_\Q(s)=\pi^{-s/2}\Gamma(s/2)\zeta(s)$ is completed. At every place we use the adjoint convention
$L_v(s,\Pi_v,\Ad)=L_v(s,(\Sym^3\pi_f)_v,\Sym^2)$
of \cite[(1.20)]{FurusawaMorimoto2021}. The standing hypotheses make the denominator nonzero, so $\mathscr L_d(\Pi)\geq0$.

The normalized distribution estimate in Theorem~\ref{thm:pointwise-distribution} follows from the unconditional polynomial bound below and the arithmetic estimates of Proposition~\ref{prop:arithmetic-moment}.

\begin{theorem}\label{thm:L-bound}
With the above notations, we have
\begin{equation}\label{eq:L-bound}
 \frac{|\mathcal B^\natural_{d,b}(\varphi)|^2}
 {\langle\varphi,\varphi\rangle}
 \ll_\varphi H(b)^{3}\mathscr L_d(\Pi)
\end{equation}
for every negative fundamental discriminant
$d\equiv a_0\pmod{M_\varphi}$ and every $b\in\Q^\times$.
\end{theorem}

\begin{remark}\label{rem:nonvanishing-scope}
For $\gcd(d,M_\Pi)=1$, $\epsilon(\tfrac12,\Pi\otimes\chi_d)
 =\epsilon(\tfrac12,\Pi)\chi_d(M_\Pi)$.
Thus the root number is constant on the chosen progression;
if it is negative, then $L(1/2,\Pi\otimes\chi_d)=0$.

If $f$ has level $1$, all finite epsilon factors are $1$, while
\cite[\S~5]{Prasad2007} gives
$\epsilon_\infty(1/2,\Pi)=i^{3\kappa-2}i^\kappa=-1$.
Hence $L(1/2,\Pi)=0$.
In either case, $\mathscr L_d(\Pi)=0$, and all Bessel periods
vanish by Theorem~\ref{thm:L-bound} and
Proposition~\ref{prop:zero-and-negative}.
\end{remark}

Section~\ref{sec:geometry} gives the CM interpretation and the qualitative
Duke limit for fixed base functions. The quantitative bounds use
the finite product and real reference factor of
Section~\ref{sec:local}. Section~\ref{sec:global} proves
Theorem~\ref{thm:L-bound}, the arithmetic coefficient bounds, and raw
exponential decay; Section~\ref{sec:Fourier} sums the coefficients to
prove Theorems~\ref{thm:distribution} and~\ref{thm:pointwise-distribution}.

\section{Bessel quotients and CM averages}
\label{sec:geometry}

\subsection{Notations and preliminaries}
\label{subsec:standing-data}

Write $\A=\R\times\A_{\mathrm{fin}}$ and $\widehat\Z=\prod_{p<\infty}\Z_p$.
Compact additive adelic quotients carry probability Haar measure unless
specified otherwise, and $\,^{\mathrm T}A$ denotes transpose.
We use right translation $R(g)\varphi(x)=\varphi(xg)$.
For $b\in F^\times$, where $F=\Q$, $\R$, or $\Q_p$, set
\[
 \delta(b)=\operatorname{diag}(bI_2,I_2).
\]
The similitude
multiplier $\mu$ and the complexified Lie algebra are given by
\[
 {}^{\mathrm T}g\mathbf Jg=\mu(g)\mathbf J,\qquad
 \mathbf J=\begin{pmatrix}0&I_2\\-I_2&0\end{pmatrix},\qquad
 \mathfrak g_\C=\operatorname{Lie}(\GSp(4,\R))\otimes_\R\C.
\]

Let $f$ be a non-CM holomorphic Hecke newform of even weight
$\kappa\geq2$ and trivial central character, and let $\pi_f$ be its associated
cuspidal automorphic representation of $\GL(2, \A)$. The transfer $\Sym^3\pi_f$ is
cuspidal by \cite[Theorem~B]{KimShahidi2002}; we fix its cuspidal
realization $\Pi$ on $\GSp(4,\A)$ with trivial central character and
holomorphic real component from
\cite[Theorem~A$'$]{RamakrishnanShahidi2007}.
Identify $K=\Sp(4,\R)\cap \mathrm O(4)$ with $\mathrm U(2)$ by
$\begin{psmallmatrix}A&B\\-B&A\end{psmallmatrix}\mapsto A+iB$.
The holomorphic constituent has lowest $K$-type of highest weight
$(2\kappa-1,\kappa+1)$,
\begin{equation}\label{e_lowest_K_type}
  W_\infty\simeq\rho_{\kappa+1,\kappa-2}
  ={\det}^{\kappa+1}\otimes\Sym^{\kappa-2}(\C^2).
\end{equation}
Extend this constituent trivially on positive scalars to
$G^+=\R_{>0}\Sp(4,\R)$, and take
\[
 \Pi_\infty\simeq B_{\kappa+1,\kappa-2}
 =\operatorname{Ind}_{G^+}^{\GSp(4,\R)}\Pi_\infty^{\mathrm{hol}},
 \qquad
 \Pi_\infty|_{G^+}\simeq
 \Pi_\infty^{\mathrm{hol}}\oplus\Pi_\infty^{\mathrm{anti}},
\]
as in \cite[\S~2.5]{Schmidt2017}.
Negative similitude exchanges the constituents. For
$\delta(-1)=\diag(-I_2,I_2)$, the type
$W_\infty^-:=\Pi_\infty(\delta(-1))W_\infty$ is antiholomorphic,
with highest weight $(-\kappa-1,1-2\kappa)$.
Thus ``holomorphic discrete series'' is understood in this full-group
sense. Since $\kappa\geq2$, $\Pi_\infty$ is tempered and non-generic.

By Deligne's theorem~\cite{Deligne1974}, the unramified Satake
parameters of $\pi_f$ and their symmetric cubes have modulus one.
We assume temperedness also at ramified finite places.
We use \cite[Theorem~1]{FurusawaMorimoto2021} with the real normalization
of Subsection~\ref{subsec:real-factor}.

Fix a factorizable smooth vector $\varphi=\otimes_v\varphi_v\in\Pi$
with $0\neq\varphi_\infty\in W_\infty$, and the standard additive character
$\psi:\Q\backslash\A\to\C^\times$, with
$\psi_\infty(x)=e^{2\pi ix}$ and $\operatorname{cond}(\psi_p)=\Z_p$.
Choose the scalar Haar measures $dx_v$ self-dual with respect to $\psi_v$, so that $\vol(\Z_p)=
 \vol(\Q\backslash\A)=1$.
Global automorphic quotient measures are Tamagawa measures unless a
probability normalization is specified. Local measures remain fixed as
the Fourier frequency varies.

All Bessel periods, models, and relative characters have trivial torus
character $\Lambda=\triv$; their additive characters are specified
separately. We use invariant Hermitian inner products, linear in the first variable.
The height $H(\cdot)$ and $L$-function conventions are those of the introduction. Unless otherwise indicated, implied constants may depend on the fixed
data but are independent of $b,d,D$.

\subsection{The Bessel period}
\label{subsec:scalar-Bessel-quotient}

Let $V=\Sym_2$ be the additive group of symmetric matrices.
For every commutative ring $A$, let
\[
 V(A)=\{X\in M_2(A):\,^{\mathrm{T}}X=X\},\qquad
 N(A)=\left\{n(X)=\begin{pmatrix}I_2&X\\0&I_2\end{pmatrix}:X\in V(A)\right\}.
\]
For nonsingular $S\in V(F)$, with $F$ of characteristic zero, define
the algebraic $F$-torus $T_S$ by its points over an $F$-algebra $A$:
\begin{equation}\label{e:Ts_defn}
   T_S(A)=\{t\in\GL(2,A):\,^{\mathrm{T}}tSt=\det(t)S\}.
\end{equation}
For a negative fundamental discriminant $d$, set
\begin{equation}\label{eq:Sd}
 S_d=
 \begin{cases}
 \begin{pmatrix}-d/4&0\\0&1\end{pmatrix},&d\equiv0\pmod{4},\\[3mm]
 \begin{pmatrix}(1-d)/4&1/2\\1/2&1\end{pmatrix},&d\equiv1\pmod{4}.
 \end{cases}
\end{equation}
Thus $S_d$ is positive definite and $4\det S_d=|d|$.
Let $T_d=T_{S_d}\simeq\operatorname{Res}_{K_d/\Q}\mathbb G_m$, and embed it in $\GSp(4)$ through
\[
 m(t)=\begin{pmatrix}t&0\\0&\det(t)\,^{\mathrm{T}}(t^{-1})\end{pmatrix}.
\]
Let $R_d=N\rtimes m(T_d)$ be the Bessel subgroup, and let $Z_G$ be the center of $\GSp(4)$.
Then
\begin{equation}\label{eq:compact_quotient}
 [R_d]_Z=Z_G(\A)R_d(\Q)\backslash R_d(\A)
\end{equation}
is compact: $x\mapsto x/\bar x$ identifies $T_d/Z_G$ with the anisotropic
norm-one torus $\operatorname{Res}^{1}_{K_d/\Q}\mathbb G_m$, and
$N(\Q)\backslash N(\A)$ is compact.
With the Tamagawa measures of
\cite{FurusawaMorimoto2021,HuaMiao2025}, put
\begin{equation*}
  \Tscr_d=Z_G(\A)T_d(\Q)\backslash T_d(\A),\qquad  v_d=\vol(\Tscr_d)=2.
\end{equation*}
The volume $2$ is the Tamagawa number of $T_d/Z_G$ \cite[\S~5, Corollary]{Ono1963}.
Write 
\begin{equation}
d\nu_{T,d}=dt/2,\qquad d\nu_{R,d}=dr/2   
\end{equation}
for probability Haar measures.
On the unipotent subgroup, use the coordinate measures
\begin{equation}\label{eq:unipotent-coordinate-measure}
 X=\begin{pmatrix}x_{11}&x_{12}\\x_{12}&x_{22}\end{pmatrix},
 \qquad dX_v=dx_{11,v}\,dx_{12,v}\,dx_{22,v}.
\end{equation}
Here $dX_\infty$ is coordinate Lebesgue measure, $\vol(N(\Z_p))=1$,
and $dX=\prod_vdX_v$ is Tamagawa measure with quotient volume $1$
\cite[\S~1.1]{FurusawaMorimoto2021}.
For the trace pairing, the self-dual measure is $|2|_v^{1/2}dX_v$.
At each finite prime $p$, normalize $dt_p$ on
$\Q_p^\times\backslash T_d(\Q_p)$ so that the image of
$T_d(\Q_p)\cap\GL(2,\Z_p)$ has volume $1$.
For $S\in\Sym_2(\R)$ positive definite, let $dt_S$ be probability
Haar measure on $\R^\times\backslash T_S(\R)$.
Write $dt_\infty$ for this measure when $S=S_d$.
Conjugation by $m(t)$ acts on $V$ as
$\rho(t)X=(\det t)^{-1}tX\,^{\mathrm{T}}t$. By \eqref{e:Ts_defn} we see that
\begin{equation}\label{eq:trace-invariance}
 \tr(S_d\rho(t)X)=\tr(S_dX).
\end{equation}
The action $X\mapsto tX\,^{\mathrm{T}}t$ on $\Sym_2$ has determinant
$(\det t)^3$, so
\[\det\rho(t)=(\det t)^{-3}(\det t)^3=1.\]
Thus the group law
$(X,t)(Y,s)=(X+\rho(t)Y,ts)$
preserves $dX\,dt$ on both sides. Hence $R_d$ is unimodular, with quotient
measure disintegrating over $\Tscr_d$ and $N(\Q)\backslash N(\A)$.
For $b\in\Q$, the character
$\theta_{d,b}(n(X))=\psi\bigl(b\tr(S_dX)\bigr)$
therefore extends to
\begin{equation*}
  \Theta_{d,b}(n(X)m(t))=\theta_{d,b}(n(X))
\end{equation*}
on $R_d(\A)$,
with trivial torus character, and is trivial on $Z_G(\A)R_d(\Q)$.
Define the Bessel period
\begin{equation}\label{eq:Bessel-period}
 \mathcal B_{d,b}(\varphi)=
 \int_{[R_d]_Z}\varphi(r)\Theta_{d,b}(r)^{-1}\,dr.
\end{equation}

\begin{proposition}
\label{prop:zero-and-negative}
For every negative fundamental discriminant $d$ and $b\leq0$,
$\mathcal B_{d,b}(\varphi)=0$.
\end{proposition}

\begin{proof}
For $b=0$, the integral \eqref{eq:Bessel-period} becomes
\[
 \mathcal B_{d,0}(\varphi)
 =\int_{\Tscr_d}\int_{N(\Q)\backslash N(\A)}
 \varphi(n(X)m(t))\,dX\,dt.
\]
The inner integral is the constant term of $\varphi$ along the
Siegel unipotent subgroup and vanishes by cuspidality.
For $b<0$, the matrix $bS_d$ is negative definite.  With $\psi_\infty(x)=e^{2\pi ix}$, the choice
$\varphi_\infty\in W_\infty\subset\Pi_\infty^{\mathrm{hol}}$ restricts
its holomorphic Fourier indices to positive-semidefinite matrices,
and cuspidality excludes singular indices.  Thus the Fourier coefficient
indexed by $bS_d$, and hence its torus integral, vanishes.
\end{proof}

Put $\mathscr X=\Q\backslash\A$.  The scalar map
$\iota_d:[R_d]_Z\to\mathscr X$ is defined by
\begin{equation}\label{e_iota_d}
  \iota_d([n(X)m(t)])=\tr(S_dX)\pmod{\Q}.
\end{equation}
This is well-defined by \eqref{eq:trace-invariance} and the
rationality of $S_d$. Define the finite complex Borel measure on $\mathscr X$ by
\begin{equation}
  \mathscr A_d(\varphi)
 =(\iota_d)_*\bigl(\varphi(r)\,d\nu_{R,d}(r)\bigr)
 =(\iota_d)_*\bigl(\varphi(r)\,dr/2\bigr).
\end{equation}
Its Fourier coefficients are
\begin{equation}\label{eq:CM-average-Fourier}
 \mathcal C_{d,b}(\varphi)=\int_{\mathscr X}\psi(-bx)\,d\mathscr A_d(\varphi)(x)
 =\frac{\mathcal B_{d,b}(\varphi)}{2}.
\end{equation}
For a finite complex Borel measure $\mu$ on $\mathscr X$,
its total variation norm is
\begin{equation*}
   \|\mu\|_{\mathrm{TV}}
 :=\sup_{\substack{
       \Phi:\mathscr X\to\C\ \text{bounded Borel}\\
       \|\Phi\|_\infty\le1}}
 \left|\int_{\mathscr X}\Phi\,d\mu\right|.
\end{equation*}
Compactness gives $\|\mathscr A_d(\varphi)\|_{\mathrm{TV}}
\leq\int_{[R_d]_Z}|\varphi|\,d\nu_{R,d}$, and
Proposition~\ref{prop:zero-and-negative} gives total mass zero.

\begin{lemma}\label{lem:scalar-ray}
Let $S,R\in\Sym_2(\Q)$ with $S$ nonsingular.  The character
$n(X)\mapsto\psi(\tr(RX))$ on $N(\A)$ is invariant under $T_S(\A)$
if and only if $R\in\Q S$.
\end{lemma}

\begin{proof}
If $R\in\Q S$, the required invariance follows from
\eqref{eq:trace-invariance}, with $S_d$ replaced by $S$.
Conversely, character invariance and the nondegenerate adelic additive
pairing imply that $X\mapsto\tr(RX)$ is fixed by $T_S(\Q)$, hence by
$T_S$ because $T_S(\Q)$ is Zariski dense.
Over $\overline\Q$, take $S$ off-diagonal and $T_S$ diagonal.
The contragredient action of $t=\diag(a,b)$ on $V^\vee$ has weights
$b/a,1,a/b$, so the invariant subspace is one-dimensional.
By \eqref{e:Ts_defn}, it contains $X\mapsto\tr(SX)$; descent therefore gives
$(V^\vee)^{T_S}=\Q\tr(S\,\cdot)$.
Hence $R\in\Q S$. This completes the proof.
\end{proof}

\subsection{The Duke base and the Bessel fiber}
\label{subsec:Duke-base}

Inside the Siegel parabolic, let
$M^1=\{m(t):t\in\GL(2)\}$ and $P^1=NM^1$.
The map $t\mapsto m(t)$ identifies $M^1$ with $\GL(2)$ and
$Z_G\backslash M^1$ with $\PGL(2)$; the Levi projection
is $\operatorname{pr}_M(n(X)m(t))=m(t)$.
Set $K_\infty^M=m(\SO(2))$ and
$K_{\mathrm{fin}}^M=m(\GL(2,\widehat\Z))$.  The adelic modular quotient is
$$\mathcal Y_M^{\mathrm{ad}}= Z_G(\A)M^1(\Q)\backslash M^1(\A)/(K_\infty^MK_{\mathrm{fin}}^M).$$
It is represented by
$\mathcal Y_M=\{[m(g_\infty),1_{\mathrm{fin}}]:\det g_\infty>0\}$.
By strong approximation for $\SL(2)$~\cite[Theorem~7.12]{PlatonovRapinchuk1994},
for every $g_{\mathrm{fin}}\in\GL(2,\A_{\mathrm{fin}})$ there exists
$\gamma\in\GL(2,\Q)$ such that
\begin{equation}\label{eq:finite-component-reduction}
 \det\gamma>0,\qquad
 \gamma g_{\mathrm{fin}}\in\GL(2,\widehat\Z).
\end{equation}
Since $\diag(-1,1)\in\GL(2,\Z)$, representatives may be chosen
with positive real determinant. The remaining rational identifications
are induced by $\SL(2,\Z)$, giving
\[
 \SL(2,\Z)\backslash\mathbb H\simeq\mathcal Y_M,\qquad
 [g_\infty\cdot i]\longmapsto[m(g_\infty),1_{\mathrm{fin}}].
\]
Through this identification, we equip $\mathcal Y_M$ with the
hyperbolic probability measure $d\mu_M$.

Write $S_d=\left(\begin{smallmatrix}a_d&b_d/2\\b_d/2&c_d\end{smallmatrix}\right)$,
and set
\[
 x_d=-\frac{b_d}{2a_d},\qquad y_d=\frac{\sqrt{|d|}}{2a_d},\qquad
 z_d=x_d+iy_d,\qquad
 g_{d,\infty}=\begin{pmatrix}\sqrt{y_d}&x_d/\sqrt{y_d}\\0&1/\sqrt{y_d}\end{pmatrix}.
\]
Then $g_{d,\infty}\cdot i=z_d$ and
$S_d=\sqrt{\det S_d}\,\,^{\mathrm{T}}(g_{d,\infty}^{-1})g_{d,\infty}^{-1}$.
For $h=g_{d,\infty}^{-1}tg_{d,\infty}$, $t\in T_d(\R)$, this gives
$\,^{\mathrm T}hh=\det(h)I_2$ with $\det h>0$, hence
\begin{equation}\label{eq:compact-real-torus}
 g_{d,\infty}^{-1}T_d(\R)g_{d,\infty}
 \subset\R_{>0}\SO(2).
\end{equation}
Let $\mathbf g_d=(m(g_{d,\infty}),1_{\mathrm{fin}})$, and define
\begin{equation}\label{e:eta_d}
  \eta_d:\Tscr_d\longrightarrow\mathcal Y_M,\quad [t]\longmapsto[m(t)\mathbf g_d].
\end{equation}
\begin{lemma}
\label{lem:CM-packet}
Let $\mathcal O_{K_d}$ be the ring of integers of $K_d$. The map $\eta_d$ defined in \eqref{e:eta_d} is well-defined, and its image is the packet
$\operatorname{CM}(d)$ indexed by the proper ideal classes of
$\mathcal O_{K_d}$.  Its probability Haar pushforward is
\begin{equation}\label{eq:CM-packet-measure}
 (\eta_d)_*\nu_{T,d}=\mu_d
 =\frac1{h(K_d)}\sum_{z\in\operatorname{CM}(d)}\delta_z,
\end{equation}
where $h(K_d)$ is the class number and $\delta_z$ denotes the unit point mass at $z$.
\end{lemma}

\begin{proof}
The left quotient absorbs $T_d(\Q)$ and $Z_G(\A)$.
By \eqref{eq:compact-real-torus}, the map $\eta_d$ factors through
\[
 T_d(\Q)\backslash T_d(\A_{\mathrm{fin}})/
 \bigl(\A_{\mathrm{fin}}^\times K_{T,d}\bigr),\qquad
 K_{T,d}=T_d(\A_{\mathrm{fin}})\cap\GL(2,\widehat\Z).
\]
The matrix $S_d$ is $\SL_2(\Z)$-equivalent to the matrix $S_E$ in
\cite[(8.2.11)]{FurusawaMorimoto2024}:
$S_d={}^{\mathrm T}u_dS_Eu_d$ for
$u_d=\begin{psmallmatrix}-b_d&-1\\1&0\end{psmallmatrix}$.
The standard maximal-order embedding therefore identifies $K_{T,d}$
with $\widehat{\mathcal O}_{K_d}^{\times}$.
Since $\A_{\mathrm{fin}}^\times=\Q^\times\widehat\Z^\times$, the above quotient is
\[
 K_d^\times\backslash\A_{K_d,\mathrm{fin}}^\times/
 \widehat{\mathcal O}_{K_d}^{\times}
 \simeq\Cl(\mathcal O_{K_d}).
\]
Here $\A_{K_d,\mathrm{fin}}$ denotes the finite adele ring of $K_d$.
As $g_{d,\infty}\cdot i=z_d$, the usual proper ideal-class
correspondence identifies its image under $\eta_d$ with the full
packet $\operatorname{CM}(d)$; see
\cite[\S~8.2.3]{FurusawaMorimoto2024}.
Finally, the fibers of $\eta_d$ are the $h(K_d)$ cosets of the image
of $T_d(\R)K_{T,d}$ in $\Tscr_d$.
Translation invariance and total mass one assign mass $1/h(K_d)$
to each fiber, proving \eqref{eq:CM-packet-measure}.
\end{proof}

\begin{theorem}[\text{Duke's equidistribution theorem \cite{Duke1988}}]
\label{thm:Duke-base}
As $|d|\to\infty$ through negative fundamental discriminants,
$\mu_d$ converges weakly to $\mu_M$.  In particular, for every
$F\in C_c^\infty(\mathcal Y_M)$,
\begin{equation*}
 \frac1{h(K_d)}\sum_{z\in\operatorname{CM}(d)}F(z)
 \longrightarrow\int_{\mathcal Y_M}F\,d\mu_M.
\end{equation*}
The same limit holds along any fixed congruence progression of such discriminants.
\end{theorem}
The projection
\[p_d:[R_d]_Z\longrightarrow\Tscr_d,\quad [n(X)m(t)]\longmapsto[t]\]
is well-defined by the rational and central quotients. Its fibers
identify with $V(\Q)\backslash V(\A)$.

\begin{proposition}
\label{prop:cylindrical-Duke-lift}
Let $\beta_d=\eta_d\circ p_d$. Then, with probability Haar
measure $dx$ on $\mathscr X$,
\begin{equation*}
 (\beta_d,\iota_d)_*\nu_{R,d}=\mu_d\otimes dx.
\end{equation*}
Consequently, for $F\in C_c^\infty(\mathcal Y_M)$ and
$\Phi\in C(\mathscr X)$, as $|d|\to\infty$,
\begin{equation*}
 \int_{[R_d]_Z}F(\beta_d(r))\Phi(\iota_d(r))\,d\nu_{R,d}(r)
 \longrightarrow
 \left(\int_{\mathcal Y_M}F\,d\mu_M\right)
 \left(\int_{\mathscr X}\Phi\,dx\right).
\end{equation*}
\end{proposition}

\begin{proof}
Unimodularity gives Weil's integration formula
\[
 \int_{[R_d]_Z}H(r)\,d\nu_{R,d}(r)
 =\int_{\Tscr_d}\int_{V(\Q)\backslash V(\A)}
 H(n(X)m(t))\,dX\,d\nu_{T,d}(t).
\]
The inner probability Haar measure is independent of the representative
of $[t]$, and the adelic center acts trivially.
The trace map
$X\mapsto\tr(S_dX)$ is a surjective homomorphism onto $\mathscr X$
because $(S_d)_{22}=1$; its Haar pushforward is $dx$.
For $H(r)=F(\beta_d(r))\Phi(\iota_d(r))$, the inner integral is
\[
 F(\eta_d([t]))\int_{\mathscr X}\Phi(x)\,dx.
\]
Lemma~\ref{lem:CM-packet} identifies the remaining torus integral with
$\int F\,d\mu_d$.  Finite linear combinations of product functions
determine Radon measures on $\mathcal Y_M\times\mathscr X$, giving
the product identity. The limit follows immediately from Theorem~\ref{thm:Duke-base}.
\end{proof}

To express the Fourier coefficients \eqref{eq:CM-average-Fourier} of $\mathscr A_d(\varphi)$ as a class group average, define
\[
 \mathcal W_{d,b}\varphi([t])=
 \int_{N(\Q)\backslash N(\A)}
 \varphi(n(X)m(t))
 \psi\bigl(-b\tr(S_dX)\bigr)\,dX.
\]
Automorphy, \eqref{eq:trace-invariance} and the trivial central
character make this a well-defined function on $\Tscr_d$.

Each fiber $\eta_d^{-1}(z)$ has measure $1/h(K_d)$; let
\[\kappa_{d,z}=h(K_d)\nu_{T,d}|_{\eta_d^{-1}(z)}\]
be its conditional
probability measure. Disintegrating the Bessel period over these
fibers and using $dt=2\,d\nu_{T,d}$ gives
\begin{equation}\label{eq:Bessel-as-CM-average}
 \begin{aligned}
 \mathcal C_{d,b}(\varphi)
 =\frac{\mathcal B_{d,b}(\varphi)}{2}
 &=\int_{\Tscr_d}\mathcal W_{d,b}\varphi(t)\,d\nu_{T,d}(t)\\
 &
 =\frac1{h(K_d)}\sum_{z\in\operatorname{CM}(d)}
 \int_{\eta_d^{-1}(z)}
 \mathcal W_{d,b}\varphi(t)\,d\kappa_{d,z}(t).
 \end{aligned}
\end{equation}
Thus the period averages the $d,b$-dependent observable
$\mathcal W_{d,b}\varphi$ over the conditional fibers above the CM packet.

\section{Local Bessel factors and polynomial control}
\label{sec:local}

For the finite-place arguments in this section, let $p$ be a prime and $F=\Q_p$, with ring of integers $\OF=\Z_p$, maximal ideal $\mathfrak p=p\Z_p$, uniformizer $\varpi=p$, residue cardinality $q=p$ and $|p|_F=p^{-1}$. Take $\psi_F=\psi_p$ and $dX=dX_p$ from
\eqref{eq:unipotent-coordinate-measure}.

\subsection{Finite local factors and uniformity}

The trace pairing on $V=\Sym_2$ is non-degenerate, also when $p=2$.
Let $\Pi_F$ be an irreducible unitarizable admissible representation of
$\GSp(4,F)$.
Fix a nonsingular $S\in V(F)$. Recall that the torus character is trivial, and write
$$\theta_S(n(X))=\psi_F(\tr(SX)),\quad  R_S=N\rtimes m(T_S),\quad  \chi_S(n(X)m(t))=\theta_S(n(X)).$$
We identify $T_S$ with its image under $m$ when convenient.
A special Bessel form on $\Pi_F$ is a sesquilinear form
$\mathcal H_S$ satisfying
\begin{equation}
 \label{eq:Bessel-form-equivariance}
 \mathcal H_S(\Pi_F(r_1)v,\Pi_F(r_2)w)
 =\chi_S(r_1)\overline{\chi_S(r_2)}\mathcal H_S(v,w)
 \qquad(r_1,r_2\in R_S(F)).
\end{equation}
For $b\in F^\times$, define
\begin{equation}
 \label{eq:transported-form}
 \mathcal I_{bS}(v,w)
 =|b|_F^{-3}\mathcal H_S
 \bigl(\Pi_F(\delta(b))v,\Pi_F(\delta(b))w\bigr).
\end{equation}
Here the Haar measure on $V(F)$ is kept fixed as $b$ varies.

The stable integral $\int_{V(F)}^{\mathrm{st}}f(X)\,dX$ is the common
value of $\int_Uf(X)\,dX$ over compact open subgroups $U\supset U_0$
for some compact open $U_0$; see \cite[\S~1.3.1, Definition~2]{FurusawaMorimoto2021}.
Here $U_0$ may depend on the vectors and torus element.
For a nonzero smooth vector $\phi$, $b\in F^\times$,
$t\in T_S(F)$, and an integer $k\geq0$, set
\[
 \Phi_\phi(g)=\frac{\langle\Pi_F(g)\phi,\phi\rangle}{\|\phi\|^2},
 \qquad
 J_{0,b}^{(k)}(\phi,t)=
 \int_{\mathfrak p^{-k}V(\OF)}
 \Phi_\phi(n(X)t)\psi_F(-b\tr(SX))\,dX.
\]
Whenever the stable integral exists, the definition
\begin{equation*}
 J_{0,b}(\phi,t)=\lim_{k\to\infty}J_{0,b}^{(k)}(\phi,t),
 \qquad J_0=J_{0,1},
\end{equation*}
is meaningful. 

\begin{lemma}
\label{lem:measure-covariance}
If $b=u\varpi^j$, $u\in\OF^\times$, and $j\geq0$, then
\begin{align}
 J_{0,b}^{(k)}(\phi,t)
 &=q^{3j}J_{0,1}^{(k-j)}(\Pi_F(\delta(b))\phi,t)
 \qquad(k>j),\label{e:finite-covariance-psiM}\\
 J_{0,b}(\phi,t)
 &=q^{3j}J_0(\Pi_F(\delta(b))\phi,t).
 \label{e:stable-covariance-psiM}
\end{align}
\end{lemma}

\begin{proof}
For $b\in F^\times$, we have
$\delta(b)n(X)\delta(b)^{-1}=n(bX)$, and $\delta(b)$ commutes
with every element of $T_S(F)$. Unitarity therefore gives
\[
 \Phi_\phi(n(b^{-1}Y)t)
 =\Phi_{\Pi_F(\delta(b))\phi}(n(Y)t).
\]
The substitution $Y=bX$ has Jacobian $dX=|b|_F^{-3}\,dY$.
For $b=u\varpi^j$, it sends $\mathfrak p^{-k}V(\OF)$ to
$\mathfrak p^{-(k-j)}V(\OF)$, proving
\eqref{e:finite-covariance-psiM}.
Multiplication by $b$ also maps the compact open subgroups
containing $U_0$ onto those containing $bU_0$.
Passing to stable integrals proves
\eqref{e:stable-covariance-psiM}.
\end{proof}

\begin{proposition}
\label{prop:fixed-place}
For every finite-dimensional space $E$ of smooth vectors of $\Pi_F$,
there is a constant $C_E$ such that
\begin{equation}
 \label{eq:fixed-place-bound}
 |\mathcal I_{bS}(v,w)|
 \leq C_E\|v\|\|w\|\max\{1,|b|_F^{-1}\}^{3}
 \qquad(v,w\in E,\ b\in F^\times).
\end{equation}
\end{proposition}

\begin{proof}
Choose $r\geq1$ such that $E\subset E_r:=\Pi_F^{K_r}$ and
$\theta_S$ is trivial on $n(\mathfrak p^rV(\OF))$, where
\[
K_r=\ker\bigl(\GSp(4,\OF)\longrightarrow \GSp(4,\OF/\mathfrak p^r)\bigr).
\]
By admissibility, $E_r$ is finite-dimensional, so the restriction
of $\mathcal H_S$ to $E_r\times E_r$ is bounded:
\[
 |\mathcal H_S(x,y)|\leq C_r\|x\|\|y\|
 \qquad(x,y\in E_r).
\]
Put $h=\delta(\varpi)$. Using probability Haar measure $dk$ on $K_r$,
define
\[
 P_r=\int_{K_r}\Pi_F(k)\,dk.
\]
Then $P_r$ is the orthogonal projection onto $E_r$. Write $b=u\varpi^j$ with $u\in\OF^\times$ and $j\in\Z$.
Since $\delta(u)$ normalizes $K_r$ and commutes with $h$,
$\Pi_F(\delta(u))$ acts unitarily on $E_r$.
It therefore suffices to consider $b=\varpi^j$.

For $j\geq0$, one has
$$K_rh^jK_r =\coprod_{X\in\mathcal R_j}n(X)h^jK_r,$$
where $\mathcal R_j$ represents
$\mathfrak p^rV(\OF)/\mathfrak p^{r+j}V(\OF)$. Evidently, $|\mathcal R_j|=q^{3j}$, and so
\[
 P_r\Pi_F(h^j)v=q^{-3j}\sum_{X\in\mathcal R_j}\Pi_F(n(X)h^j)v.
\]
The character is trivial on these representatives, so equivariance in
both variables, applied to the finite sums, gives
\begin{align*}
 |\mathcal H_S(\Pi_F(h^j)v,\Pi_F(h^j)w)|
 &=|\mathcal H_S(P_r\Pi_F(h^j)v,P_r\Pi_F(h^j)w)|\\
 &\leq C_r\|P_r\Pi_F(h^j)v\|\|P_r\Pi_F(h^j)w\|
 \leq C_r\|v\|\|w\|.
\end{align*}

For every integer $j$, the vector $\Pi_F(h^j)v$ is fixed by
$n(\mathfrak p^{r+j}V(\OF))$. For all sufficiently negative $j$,
$\theta_S$ is nontrivial on this subgroup, so equivariance forces
\[\mathcal H_S(\Pi_F(h^j)v,\Pi_F(h^j)w)=0.\]
The finitely many remaining negative $j$ give bounded forms on
$\Pi_F(h^j)E_r$, so the same bound holds for all $j\in\Z$ after
increasing $C_r$. Applying it to $\Pi_F(\delta(u))v$ and
$\Pi_F(\delta(u))w$, which lie in $E_r$ with unchanged norms,
\eqref{eq:transported-form} gives
$|\mathcal I_{bS}(v,w)|\leq C_E|b|_F^{-3}\|v\|\|w\|$
and hence \eqref{eq:fixed-place-bound}.
\end{proof}
We now assume that $\Pi_F$ has trivial central character.
For smooth vectors $x,y$, define
\begin{equation}
 \label{eq:local-relative-integral}
 \alpha_{S,b}(x,y)
 =\int_{F^\times\backslash T_S(F)}\int_{V(F)}^{\mathrm{st}}
 \langle\Pi_F(m(t)n(X))x,y\rangle
 \psi_F(-b\tr(SX))\,dX\,dt.
\end{equation}
If $\Pi_F$ is tempered, the inner integral in
\eqref{eq:local-relative-integral} stabilizes for every
$t\in T_S(F)$, and the outer torus integral converges absolutely;
see \cite[Proposition~3.1 and Theorem~2.1(1)]{LiuYifeng2016}.
For general unitary data, we assume these two properties.
Since $T_{bS}=T_S$ by \eqref{e:Ts_defn}, this form satisfies
\eqref{eq:Bessel-form-equivariance} with $\chi_S$ replaced by
$\chi_{bS}(n(X)m(t))=\psi_F(b\tr(SX))$.
Let $K$ be the quadratic algebra determined by $S$ and $\chi_{K/F}$
its quadratic character, trivial in the split case. The torus character
is trivial, so $\AI(\triv)=\triv\boxplus\chi_{K/F}$. Set
\begin{align*}
 J_{0,b}(\phi)
 &=\int_{F^\times\backslash T_S(F)}
 J_{0,b}(\phi,t)\,dt
 =\frac{\alpha_{S,b}(\phi,\phi)}
 {\|\phi\|^2},\\
 J_b(\phi)
 &=\frac{L(1,\Pi_F,\Ad)L(1,\chi_{K/F})}
 {\zeta_F(2)\zeta_F(4)L(1/2,\Pi_F)
 L(1/2,\Pi_F\otimes\chi_{K/F})}
 J_{0,b}(\phi),
\end{align*}
assuming that the normalizing factors are finite, and write $J=J_1$.

For $S=S_d$, abbreviate \eqref{eq:local-relative-integral} by
$\alpha_{p,d,b}$. For $x\neq0$, the normalized local factor is
\begin{equation}
 \label{eq:normalized-local-factor}
 J_{p,d,b}(x)=\mathcal N_{p,d}
 \frac{\alpha_{p,d,b}(x,x)}{\langle x,x\rangle_p},
 \qquad
 \mathcal N_{p,d}=
 \frac{L(1,\Pi_p,\Ad)L(1,\chi_{d,p})}
 {\zeta_p(2)\zeta_p(4)L(1/2,\Pi_p)
 L(1/2,\Pi_p\otimes\chi_{d,p})}.
\end{equation}

Under the exceptional isomorphism $\PGSp(4)\simeq\SO(5)$, $J_{p,d,b}(x)$ equals
$\alpha_p^\natural(x,x)/\langle x,x\rangle_p$ in
\cite[Definition~5 and (1.14)]{FurusawaMorimoto2021},
with $n=2$, Bessel matrix $bS_d$, and our fixed Haar measures.

For $b\ne0$, substitute $Y=bX$ in
\eqref{eq:local-relative-integral}, so that
$dX=|b|_p^{-3}\,dY$. Unitarity gives
\begin{equation}
 \label{eq:cross-relative-transport}
 \alpha_{p,d,b}(x,y)
 =|b|_p^{-3}\,\alpha_{p,d,1}
 \bigl(\Pi_p(\delta(b))x,\Pi_p(\delta(b))y\bigr).
\end{equation}

We first bound the denominators of the frequencies $b$ for which
$\mathcal B_{d,b}\ne0$, uniformly in $d$.
\begin{lemma}
\label{lem:denominator}
Choose principal congruence depths $r_p\ge0$ fixing $\varphi_p$,
with $r_p=0$ outside $\Sigma_\varphi$, and put
\[
 Q_\varphi=\prod_{p\in\Sigma_\varphi}p^{r_p}.
\]
Then, for every negative fundamental discriminant $d$,
\[
 \mathcal B_{d,m/n}(\varphi)\neq0
 \quad\Longrightarrow\quad n\mid Q_\varphi
 \qquad(\gcd(m,n)=1,\ n>0).
\]
\end{lemma}

\begin{proof}
Embed $s_p\in R_d(\Q_p)$ adelically with identity components away
from $p$. Right invariance of the quotient measure gives
\[
 \mathcal B_{d,b}\bigl(R(s_p)\varphi\bigr)
 =\Theta_{d,b,p}(s_p)\mathcal B_{d,b}(\varphi).
\]
By the choice of $r_p$, the element $n(p^{r_p}xE_{22})$ fixes
$\varphi_p$ for every $x\in\Z_p$, where $E_{22}=\diag(0,1)$.
Since $(S_d)_{22}=1$, if $\mathcal B_{d,b}(\varphi)\neq0$,
equivariance gives
\[
 \psi_p(bp^{r_p}x)=1\qquad(x\in\Z_p).
\]
Since $\psi_p$ has conductor $\Z_p$, this implies
$bp^{r_p}\in\Z_p$.
Writing $b=m/n$ in lowest terms with $n>0$, we obtain
$v_p(n)\le r_p$ for every prime $p$, hence $n\mid Q_\varphi$.
\end{proof}

We next bound the local Bessel factors at the primes in $\Sigma_\varphi$,
uniformly in $d$ within the fixed congruence class.
\begin{lemma}
\label{lem:d-uniformity}
For each $p\in\Sigma_\varphi$, there is a constant $C_{p,\varphi}$ such that
\begin{equation}
 \label{eq:bad-place-uniform}
 |J_{p,d,b}(\varphi_p)|
 \leq C_{p,\varphi}\max\{1,|b|_p^{-1}\}^{3}
\end{equation}
for all negative fundamental $d\equiv a_0\pmod{M_\varphi}$ and
$b\in\Q_p^\times$.
\end{lemma}

\begin{proof}
Fix $p\in\Sigma_\varphi$ and a negative fundamental discriminant
$d_\nu\equiv a_0\pmod{M_\varphi}$, and put $S_{p,\nu}=S_{d_\nu}$.
For $g\in\GL(2,\Z_p)$, write
$S^g=(\det g)^{-1}\,^{\mathrm T}gSg$.
Since $d\equiv d_\nu\equiv a_0\pmod{M_\varphi}$ and
$\gcd(a_0,M_\varphi)=1$, we have
\[
 \frac{d}{d_\nu}\in
 \begin{cases}
  1+p\Z_p,&p\ne2,\\
  1+8\Z_2,&p=2.
 \end{cases}
\]
These congruences imply that $d/d_\nu\in (\Z_p^\times)^2$. Choose $v\in\Z_p^\times$ with $v^2=d/d_\nu$ and set
\[
 g=\begin{pmatrix}v&0\\(1-v)/2&1\end{pmatrix},
 \qquad \det g=v.
\]
Then $g\in\GL(2,\Z_p)$, since $2$ is a unit for odd $p$
and $v$ is odd for $p=2$.
As $d\equiv d_\nu\equiv1\pmod4$, \eqref{eq:Sd} gives
\[
 \,^{\mathrm T}gS_{p,\nu}g
 =\begin{pmatrix}
    (1-v^2d_\nu)/4&1/2\\
    1/2&1
   \end{pmatrix}
 =S_d.
\]
Thus $S_d=vS_{p,\nu}^g$.
Since $d/d_\nu=v^2$, we also have
$\chi_{d,p}=\chi_{d_\nu,p}$ and
$\mathcal N_{p,d}=\mathcal N_{p,d_\nu}$.

Conjugation by $m(g)$ preserves the chosen measures and gives
\[
 J_{p,d,b}(\varphi_p)
 =J_{p,\nu,bv}\bigl(\Pi_p(m(g))\varphi_p\bigr),
\]
where $J_{p,\nu,bv}$ denotes the local factor for $bvS_{p,\nu}$.
Let $E$ be the span of the $\GSp(4,\Z_p)$-orbit of $\varphi_p$.
This space is finite-dimensional because $\varphi_p$ is smooth.
By Proposition~\ref{prop:fixed-place}, applied on $E$ with
$S=S_{p,\nu}$ and frequency $bv$,
\[
 |J_{p,d,b}(\varphi_p)|
 \ll_{p,\varphi}\max\{1,|bv|_p^{-1}\}^{3}
 =\max\{1,|b|_p^{-1}\}^{3},
\]
which proves \eqref{eq:bad-place-uniform}.
\end{proof}

\subsection{Spherical factors and the finite product}
\label{sec:explicit-finite}

Assume that
\[
 S=\begin{pmatrix}a&b_0/2\\b_0/2&c\end{pmatrix},\qquad
 a,b_0\in\OF,\quad c\in\OF^\times,
\]
where $d=b_0^2-4ac\neq0$ generates the discriminant ideal of the
quadratic algebra $K$; in the split case $K=F\times F$.
Write $T_S(\OF)=T_S(F)\cap\GL(2,\OF)$ and normalize $dt$ so that
its image in $F^\times\backslash T_S(F)$ has volume $1$.
By type I we mean an irreducible normalized Borel principal series
$\chi_1\times\chi_2\rtimes\sigma$ with unramified inducing characters.
Every irreducible unramified tempered representation of $\GSp(4,F)$
is of this type, with unitary inducing characters; see
\cite[Proposition~3.1]{PitaleSchmidt2009Ramanujan}.
Let $\Pi_F$ be type I and $\phi$ a nonzero
$\GSp(4,\OF)$-fixed vector as in \cite[\S~2.2]{DicksonPitaleSahaSchmidt2020}.

For $b=u\varpi^j$, $u\in\OF^\times$, $j\geq0$, put
\[
h(\ell,m)=\diag(\varpi^{\ell+2m},\varpi^{\ell+m},1,\varpi^m).
\]
Since $h(j,0)=\delta(\varpi^j)$ and $\delta(u)$ fixes $\phi$,
Lemma~\ref{lem:measure-covariance} gives
\begin{equation}\label{e:JM-spherical-covariance}
 J_{0,b}(\phi)
 =p^{3j}J_0(\Pi_F(h(j,0))\phi),\qquad J_b(\phi)
 =p^{3j}J(\Pi_F(h(j,0))\phi).
\end{equation}
The fundamental spherical calculation
\cite[Theorem~2.1(i)]{DicksonPitaleSahaSchmidt2020} gives
$J(\phi)=1$, including for ramified quadratic algebras.
By Bessel multiplicity one
\cite[\S~2.3]{KowalskiSahaTsimerman2012}, the sesquilinear form
$\alpha_{S,1}$ satisfies
\begin{equation}\label{e_rank-one_identity}
  \alpha_{S,1}(v,w)
 =C\,\beta(v)\overline{\beta(w)}
\end{equation}
for a nonzero Bessel functional $\beta$ and a constant $C$.
Since $J(\phi)=1$, we have $\beta(\phi)\neq0$, and we set
\[
 B_\phi(g)=\frac{\beta(\Pi_F(g)\phi)}{\beta(\phi)}.
\]
Since $\Pi_F$ is unitary, \eqref{e:JM-spherical-covariance} and the rank-one identity \eqref{e_rank-one_identity} give
\begin{equation}\label{e:JM-all-j}
 J_b(\phi)=p^{3j}|B_\phi(h(j,0))|^2.
\end{equation}
Let $\delta_1,\ldots,\delta_4$ be the unitary-normalized spin Satake
parameters. Set $\epsilon_{K/F}=1,0,-1$ in the split, ramified, and
unramified quadratic cases, respectively. For the trivial torus
character, Sugano's formula in
\cite[Theorem~2.5 and (2.3.7)]{KowalskiSahaTsimerman2012} yields
\begin{equation*}
 \sum_{j\geq0}B_\phi(h(j,0))\,Y^j
 =\frac{(1-p^{-2}Y)(1-\epsilon_{K/F}p^{-2}Y)}
 {\prod_{i=1}^4(1-\delta_ip^{-3/2}Y)}.
\end{equation*}
With $H_n=\sum_{\substack{a_1,\ldots,a_4\in\Z_{\geq0}\\a_1+\cdots+a_4=n}}\prod_i\delta_i^{a_i}$ for $n\geq0$
and $H_n=0$ for $n<0$, coefficient extraction gives
\begin{equation}\label{e:def-bj-Sugano}
 p^{3j/2}B_\phi(h(j,0))=H_j-p^{-1/2}(1+\epsilon_{K/F})H_{j-1}
                   +p^{-1}\epsilon_{K/F}H_{j-2}.
\end{equation}

If $\Pi_F$ is tempered, then $|\delta_i|=1$, so
$|H_n|\leq\binom{n+3}{3}$. Since $|\epsilon_{K/F}|\leq1$, \eqref{e:def-bj-Sugano} yields
\[
 |B_\phi(h(j,0))|\leq p^{-3j/2}\left[
 \binom{j+3}{3}+2p^{-1/2}\binom{j+2}{3}
 +p^{-1}\binom{j+1}{3}\right]
 \leq4(j+1)^3p^{-3j/2}.
\]
Here $\binom{a}{3}=0$ for $0\leq a<3$. Consequently, by \eqref{e:JM-all-j} we obtain
\begin{equation}\label{e:JM-polynomial-bound}
 J_b(\phi)\leq16(j+1)^6,
\end{equation}
uniformly in $p$ and the splitting type. The factor at $j=0$ is
exactly $1$.
Note that the growth order and uniformity in \eqref{e:JM-polynomial-bound}
also follow from \cite[Proposition~2.1]{ComtatMarzecBallesterosSaha2025}
with $\Lambda=1$ and $m=0$.

\begin{lemma}
\label{lem:good-place}
For a negative fundamental discriminant $d$, a prime
$p\notin\Sigma_\varphi$, and $j=v_p(b)\geq0$,
\begin{equation*}
 |J_{p,d,b}(\varphi_p)|\leq16(j+1)^6.
\end{equation*}
The factor equals $1$ when $j=0$, including when $p\mid d$.
\end{lemma}

\begin{proof}
Take $S=S_d$ and $\phi=\varphi_p$ in
\eqref{e:JM-polynomial-bound}, using $J_1(\phi)=1$ when $j=0$.
\end{proof}

\begin{proposition}
\label{prop:finite-product}
For every negative fundamental discriminant
$d\equiv a_0\pmod{M_\varphi}$ and every $b\in\Q^\times$ with
$\mathcal B_{d,b}(\varphi)\neq0$,
\begin{equation*}
 \prod_{p<\infty}|J_{p,d,b}(\varphi_p)|\ll_\varphi H(b)^{3}.
\end{equation*}
\end{proposition}

\begin{proof}
Write $b=m/n$ in lowest terms, $n>0$. Lemma~\ref{lem:denominator}
gives $n\mid Q_\varphi$, whose prime divisors lie in $\Sigma_\varphi$.
At good places $v_p(b)=v_p(m)\geq0$, and the factor is $1$ if
$p\nmid m$, including at primes dividing $d$. Put
\[
 m_{\rm good}=\prod_{p\notin\Sigma_\varphi}p^{v_p(|m|)},\qquad
 m_{\rm bad}=\prod_{p\in\Sigma_\varphi}p^{v_p(|m|)}.
\]
Since $16(j+1)^6\leq p^{3j}$ for $p\geq11$ and $j\geq1$,
Lemma~\ref{lem:good-place}, with the primes $3,5,7$ absorbed
into the implied constant, gives
\[
 \prod_{p\notin\Sigma_\varphi}|J_{p,d,b}(\varphi_p)|
 \ll m_{\rm good}^3.
\]
Lemma~\ref{lem:d-uniformity} bounds the remaining product by
$O_\varphi(m_{\rm bad}^3)$, since
$\prod_{p\in\Sigma_\varphi}\max\{1,|b|_p^{-1}\}=m_{\rm bad}$.
The assertion follows immediately from $m_{\rm good}m_{\rm bad}=|m|\leq H(b)$.
\end{proof}

\subsection{The real Bessel factor and its normalization}
\label{subsec:real-factor}

Retain $W_\infty\simeq\rho_{\kappa+1,\kappa-2}$ from
Subsection~\ref{subsec:standing-data}.
For $S>0$, use $dX=dX_\infty$ and $dt_S$ as in
Subsection~\ref{subsec:scalar-Bessel-quotient}. For smooth vectors $v,w$ of
$\Pi_\infty$ and $b\in\R^\times$, define
\begin{equation}\label{eq:real-relative-integral}
 \alpha_{S,b}(v,w)
 =\int_{\R^\times\backslash T_S(\R)}\int_{\Sym_2(\R)}
 \langle\Pi_\infty(m(t)n(X))v,w\rangle_\infty
 e^{-2\pi i b\tr(SX)}\,dX\,dt_S.
\end{equation}
By \cite[Proposition~3.15]{LiuYifeng2016}, this integral converges
absolutely and agrees with its Fourier regularization. Put $\zeta_\R(s)=\Gamma_\R(s)=\pi^{-s/2}\Gamma(s/2)$ and
\begin{equation}\label{eq:real-normalizing-factor}
 \mathcal N_\infty=
 \frac{L_\infty(1,\Pi,\Ad)L_\infty(1,\operatorname{sgn})}
 {\zeta_\R(2)\zeta_\R(4)L_\infty(1/2,\Pi)
 L_\infty(1/2,\Pi\otimes\operatorname{sgn})}>0.
\end{equation}
The real quadratic character is $\operatorname{sgn}$ for every positive-definite
$S$, so $\mathcal N_\infty$ is independent of $S$ and $b$.
For $v\neq0$, the normalized diagonal relative character is
\begin{equation*}
 J_{S,b}(v)=\mathcal N_\infty
 \frac{\alpha_{S,b}(v,v)}
 {\langle v,v\rangle_\infty}.
\end{equation*}
For $S=S_d$, write
\begin{equation}\label{eq:real-torus-measure-conversion}
 J_{\infty,d,b}(v):=J_{S_d,b}(v).
\end{equation}
With our fixed Haar measures, the same correspondence gives
$J_{\infty,d,b}(v)
=\alpha_\infty^\natural(v,v)/\langle v,v\rangle_\infty$;
see \cite[Definition~5 and (1.14)]{FurusawaMorimoto2021}.

\begin{lemma}
\label{lem:real-lowest-type}
For every positive-definite $S\in\Sym_2(\R)$, there is a continuous
special Bessel functional $\ell_S$ on the smooth vectors of
$\Pi_\infty$, with trivial torus character and unipotent character
$n(X)\mapsto e^{2\pi i\tr(SX)}$.
These functionals can be chosen so that
\[
 \mathcal N_\infty\,
 \alpha_{S,1}(v,w)
 =c(\det S)^{-3/2}\ell_S(v)\overline{\ell_S(w)}
 \qquad(v,w\text{ smooth}),
\]
where $c>0$ is independent of $S$.
Moreover, $\ell_S|_{W_\infty}\ne0$, so the normalized form on
$W_\infty$ is positive semidefinite of rank one.
For $b>0$ and $v\in W_\infty$,
\begin{equation}\label{eq:real-Bessel-function-scaling}
 \ell_S(\Pi_\infty(\delta(b))v)
 =b^{3\kappa/2}e^{-2\pi(b-1)\tr S}\ell_S(v).
\end{equation}
\end{lemma}

\begin{proof}
Put $\mathscr H_S^{\pm}=\mathcal N_\infty\,
\alpha_{S,\pm1}$. First, we consider the special case $S=I_2$.
Since $\kappa+1>2$, we may use the negative-spectrum cone realization of the
antiholomorphic constituent, as in
\cite[\S\S~1.3--1.6 and Lemma~2.1]{FrahmOlafssonOrsted2024}.
Its vectors are represented by functions on
$\Omega=\{Y\in\Sym_2(\R):Y>0\}$ with values in
$V_\rho=\rho_{\kappa+1,\kappa-2}^{\vee}$. The Hermitian density is $M(Y)\,dY$,
where $M(Y)$ is smooth and positive definite and
$dY=dy_{11}\,dy_{12}\,dy_{22}$, and
\[
 (\Pi_\infty(n(X))F)(Y)=e^{-2\pi i\tr(YX)}F(Y).
\]
Since $\kappa-2$ is even, choose $0\ne\xi\in V_\rho^{\SO(2)}$.
Choose $\chi\in C_c^\infty(\Omega)$ with $\chi(I_2)=1$ and
\[
 \chi(\gamma Y\,^{\mathrm T}\gamma)=\chi(Y)
 \qquad(\gamma\in\SO(2)),
\]
and set
\[
 F_0(Y)=\chi(Y)\xi\in C_c^\infty(\Omega,V_\rho).
\]
This is a smooth vector fixed by $\R^\times\backslash T_{I_2}(\R)$.
Put $h_0(Y)=\langle M(Y)F_0(Y),F_0(Y)\rangle_{V_\rho}$,
extended by zero outside $\Omega$. Then
\[
 \langle\Pi_\infty(n(X))F_0,F_0\rangle_\infty
 =\int_{\Sym_2(\R)}h_0(Y)e^{-2\pi i\tr(YX)}\,dY.
\]
The torus measure has mass $1$, so Fourier inversion gives
\[
 \mathscr H^-_{I_2}(F_0,F_0)
 =\frac{\mathcal N_\infty}{2}h_0(I_2)
 =\frac{\mathcal N_\infty}{2}
   \langle M(I_2)\xi,\xi\rangle_{V_\rho}>0.
\]
The vector $u_0=\Pi_\infty(\delta(-1))F_0$ lies in the holomorphic
constituent. Unitarity and $X\mapsto-X$ give
\[
 \mathscr H^+_{I_2}(u_0,u_0)
 =\mathscr H^-_{I_2}(F_0,F_0)>0.
\]
Thus $\ell_{I_2}=\mathscr H^+_{I_2}(\,\cdot\,,u_0)$ is a
nonzero continuous special Bessel functional. By multiplicity one
\cite[(1.9)]{FurusawaMorimoto2021},
\[
 \mathscr H^+_{I_2}(v,w)
 =c\,\ell_{I_2}(v)\overline{\ell_{I_2}(w)},\qquad c>0.
\]

Next, we consider the explicit Bessel functions. For $K$-finite $v$, continuity and unitarity give
\[
 |\ell_{I_2}(\Pi_\infty(g)v)|
 \ll\sum_i\|d\Pi_\infty(\operatorname{Ad}(g^{-1})D_i)v\|_\infty
 \ll_v\|g\|_{\rm PS}^{N}
\]
for finitely many $D_i\in U(\mathfrak g_\C)$, where
$\|g\|_{\rm PS}^2=|\mu(g)|^{-2}+\sum_{a,b}g_{ab}^2$.
The last bound follows from the polynomial dependence of the adjoint
coefficients on $g$ and $\mu(g)^{-1}$. Thus the growth condition in
\cite[\S~2.6.1]{PitaleSchmidt200903} is satisfied.
Let $v_{\mathrm{hw}}\in W_\infty$ be a highest-weight vector and put
\[
 B_0(g)=\ell_{I_2}(\Pi_\infty(g)v_{\mathrm{hw}}),\qquad
 h(\lambda,\zeta)=
 \diag(\lambda\zeta,\lambda\zeta^{-1},\zeta^{-1},\zeta).
\]
By continuity and the density of $K$-finite vectors, the Bessel model
map is nonzero. It is therefore injective on the irreducible holomorphic
Harish--Chandra module, so $B_0\not\equiv0$.
Taking $(l,l')=(2\kappa-1,\kappa+1)$ and trivial torus and central characters in
\cite[Theorems~3.4 and~3.10, (52), (60)--(61)]{PitaleSchmidt200903},
we may rescale $\ell_{I_2}$, and adjust $c>0$, so that
\[
 \begin{aligned}
 B_0(h(\lambda,\zeta))
 &=(2\zeta^2)^{(\kappa-2)/2}\lambda^{3\kappa/2}
   e^{-2\pi\lambda(\zeta^2+\zeta^{-2})},\\
 B_0(h(\lambda,\zeta)u)
 &=P_u(\zeta)\lambda^{3\kappa/2}
   e^{-2\pi\lambda(\zeta^2+\zeta^{-2})},\qquad u\in K,
 \end{aligned}
\]
where $\lambda,\zeta>0$ and
\[
 P_u(\zeta)=\sum_{a=2-\kappa}^{\kappa-2}A_a(u)\zeta^a,\qquad
 \sup_{u\in K}|A_a(u)|<\infty.
\]

Finally, we study the general positive-definite $S$.
Set $A=S^{1/2}$,
$g_S=\diag(A,{}^{\mathrm T}A^{-1})$, and
$\ell_S=\ell_{I_2}\circ\Pi_\infty(g_S)$.
The identities
\[
 g_Sn(X)g_S^{-1}=n(AX{}^{\mathrm T}A),\qquad
 g_Sm(T_S)g_S^{-1}=m(T_{I_2})
\]
imply that $\ell_S$ has the stated Bessel character.
Conjugation preserves the probability torus measure, and
$d(AX{}^{\mathrm T}A)=(\det S)^{3/2}dX$. Therefore, for smooth vectors $v,w$,
\begin{equation}\label{e_rank_one_identity}
 \begin{aligned}
 \mathscr H^+_S(v,w)
 &=(\det S)^{-3/2}
   \mathscr H^+_{I_2}
   (\Pi_\infty(g_S)v,\Pi_\infty(g_S)w)\\
 &=c(\det S)^{-3/2}
   \ell_S(v)\overline{\ell_S(w)}.
 \end{aligned}
\end{equation}
Write $A=\gamma_0\diag(a_1,a_2){}^{\mathrm T}\gamma_0$, where
$\gamma_0\in\SO(2)$ and $a_1\geq a_2>0$, and set
$\lambda=a_1a_2$, $\zeta=(a_1/a_2)^{1/2}$. Then
\[
 g_S=(\lambda^{-1/2}I_4)m(\gamma_0)h(\lambda,\zeta)m(\gamma_0)^{-1},
 \qquad \lambda(\zeta^2+\zeta^{-2})=\tr S.
\]
Triviality of the central and torus characters gives
\[
 \ell_S\bigl(\Pi_\infty(\delta(b)m(\gamma_0)u)v_{\mathrm{hw}}\bigr)
 =B_0(h(b\lambda,\zeta)u)\qquad(b>0,\ u\in K).
\]
For $b=1$ and $u=1$, the right side is nonzero by the first
explicit formula, so $\ell_S|_{W_\infty}\ne0$.
The second explicit formula gives
\[
 B_0(h(b\lambda,\zeta)u)
 =b^{3\kappa/2}e^{-2\pi(b-1)\tr S}B_0(h(\lambda,\zeta)u).
\]
Since $m(\gamma_0)\in K$, the vectors
$\Pi_\infty(m(\gamma_0)u)v_{\mathrm{hw}}$, $u\in K$, span $W_\infty$.
This proves \eqref{eq:real-Bessel-function-scaling}; the rank-one
identity above proves positivity and rank one on $W_\infty$.
\end{proof}
For $S=S_d$, let $A_{d,1}:W_\infty\to W_\infty$ be the unique
operator satisfying
\[
 \langle A_{d,1}v,w\rangle_\infty
 =\mathcal N_\infty\alpha_{S_d,1}(v,w)
 \qquad(v,w\in W_\infty).
\]
By Lemma~\ref{lem:real-lowest-type}, this operator is positive
semidefinite of rank one, with
\[
 \operatorname{Im}(A_{d,1})
 =\bigl(\ker(\ell_{S_d}|_{W_\infty})\bigr)^\perp.
\]
Choose a unit vector $q_d$ spanning this line and define
\begin{align*}
 \mathcal J_{\infty,d,1}
 &:=J_{\infty,d,1}(q_d)
 =\max_{\substack{v\in W_\infty\\\|v\|_\infty=1}}
 J_{\infty,d,1}(v)>0,\\
 \mathcal J_{\infty,d,b}
 &:=J_{\infty,d,b}(q_d)\qquad(b>0).
\end{align*}
These quantities are independent of the normalization of $\ell_{S_d}$
and the phase of $q_d$.

\begin{proposition}
\label{prop:real-factor}
For every negative fundamental discriminant $d$ and every $b>0$,
\begin{equation}\label{eq:real-factor}
 \frac{\mathcal J_{\infty,d,b}}
 {\mathcal J_{\infty,d,1}}
 =b^{3\kappa-3}e^{-4\pi(b-1)\tr(S_d)}.
\end{equation}
Moreover,
\[
 \mathcal J_{\infty,d,b}
 =\max_{\substack{v\in W_\infty\\\|v\|_\infty=1}}
 J_{\infty,d,b}(v)>0,
\]
and, for every nonzero $v\in W_\infty$,
\begin{equation}\label{eq:archimedean-comparison}
 J_{\infty,d,b}(v)
 =\mathcal J_{\infty,d,b}
 \frac{|\langle v,q_d\rangle_\infty|^2}{\|v\|_\infty^2},
 \qquad
 0\leq J_{\infty,d,b}(v)
 \leq\mathcal J_{\infty,d,b}.
\end{equation}
\end{proposition}

\begin{proof}
For $0\ne v\in W_\infty$, the positive rank-one operator $A_{d,1}$ gives
\[
 J_{\infty,d,1}(v)
 =\mathcal J_{\infty,d,1}
 \frac{|\langle v,q_d\rangle_\infty|^2}{\|v\|_\infty^2}.
\]
By \eqref{e_rank_one_identity}, we see that
\[
 J_{\infty,d,1}(w)
 =c(\det S_d)^{-3/2}
   \frac{|\ell_{S_d}(w)|^2}{\|w\|_\infty^2}
\]
for every nonzero smooth vector $w$ of $\Pi_\infty$, with $c>0$.
Changing variables $Y=bX$ in
\eqref{eq:real-relative-integral} and using unitarity together with
\eqref{eq:real-Bessel-function-scaling}, we obtain
\[
 \begin{aligned}
 J_{\infty,d,b}(v)
 =b^{-3}
   J_{\infty,d,1}(\Pi_\infty(\delta(b))v)
  &=cb^{-3}(\det S_d)^{-3/2}
   \frac{|\ell_{S_d}(\Pi_\infty(\delta(b))v)|^2}
        {\|v\|_\infty^2}\\
 &=b^{3\kappa-3}e^{-4\pi(b-1)\tr(S_d)}
   J_{\infty,d,1}(v).
 \end{aligned}
\]
Taking $v=q_d$ proves \eqref{eq:real-factor}.
Dividing the expression for $J_{\infty,d,b}(v)$ by its value at
$v=q_d$ gives \eqref{eq:archimedean-comparison}; Cauchy--Schwarz shows
that the maximum is attained at $q_d$.
\end{proof}

\begin{corollary}
\label{cor:real-uniform-bound}
For every negative fundamental discriminant $d$ and every $b>0$,
\begin{equation}\label{eq:real-absolute-bound}
 (\mathcal J_{\infty,d,b})^{1/2}
 \ll_{\Pi_\infty}
 b^{(3\kappa-3)/2}|d|^{\kappa-5/4}
 e^{-2\pi b\tr(S_d)}.
\end{equation}
In particular, uniformly in $d$ and every real $b>0$,
\begin{equation}\label{eq:real-height-bound}
 (\mathcal J_{\infty,d,b})^{1/2}
 \ll_{\Pi_\infty}1.
\end{equation}
\end{corollary}

\begin{proof}
For $S>0$ with eigenvalues $\Lambda\geq\mu>0$, put
$\lambda=(\det S)^{1/2}$ and $\zeta=(\Lambda/\mu)^{1/4}\geq1$.
In Lemma~\ref{lem:real-lowest-type}, the polynomial is
$O_{\Pi_\infty}(\zeta^{\kappa-2})$ uniformly on $K$. Choosing a basis of
$K$-translates of $v_{\mathrm{hw}}$ makes the bound uniform on the unit
sphere of $W_\infty$, and unitarity of $m(\gamma_0)$ makes it uniform in $\gamma_0\in\SO(2)$.
Hence
\[
 \sup_{\substack{v\in W_\infty\\\|v\|_\infty=1}}
 |\ell_{I_2}(\Pi_\infty(g_S)v)|
 \ll_{\Pi_\infty}
 \lambda^{3\kappa/2}\zeta^{\kappa-2}e^{-2\pi\tr S}.
\]
By \eqref{e_rank_one_identity}, we obtain
\[
 \max_{\substack{v\in W_\infty\\\|v\|_\infty=1}}
 J_{S,1}(v)
 \ll_{\Pi_\infty}
 (\det S)^{(3\kappa-3)/2}(\Lambda/\mu)^{(\kappa-2)/2} e^{-4\pi\tr S}.
\]
For $S_d$ in \eqref{eq:Sd}, with $D=|d|$,
\[
 \det S_d=D/4,\qquad
 D/4\leq\tr S_d\leq(D+5)/4,\qquad
 \mu\asymp1,\quad\Lambda\asymp D.
\]
Applying the bound to $bS_d$ and using the maximum identity in
Proposition~\ref{prop:real-factor}, we obtain
\[
 \mathcal J_{\infty,d,b}
 \ll_{\Pi_\infty}
 b^{3\kappa-3}D^{2\kappa-5/2}e^{-4\pi b\tr S_d}.
\]
Here $\det(bS_d)=b^2D/4$ and the eigenvalue ratio is unchanged.
Taking square roots gives \eqref{eq:real-absolute-bound}.
Finally, set $x=bD$.  Since $D\geq3$ and $\kappa\geq2$,
\[
 b^{(3\kappa-3)/2}D^{\kappa-5/4}e^{-2\pi b\tr S_d}
 \leq D^{1/4-\kappa/2}
 x^{(3\kappa-3)/2}e^{-\pi x/2}\ll_\kappa1,
\]
which proves \eqref{eq:real-height-bound} uniformly for $D\geq3$ and
$b>0$.
\end{proof}

\section{The global period and its bound}
\label{sec:global}
\subsection{Proof of the pointwise period bound}
Before proving Theorem~\ref{thm:L-bound}, we first transfer the scalar frequency to the test vector.

\begin{lemma}
\label{lem:global-covariance}
For every negative
fundamental discriminant $d$ and every $b\in\Q^\times$,
\begin{equation}
 \label{eq:global-covariance}
 \mathcal B_{d,b}(\varphi)=
 \mathcal B_{d,1}\bigl(R(\delta(b))\varphi\bigr).
\end{equation}
\end{lemma}

\begin{proof}
By automorphy,
\[
 (R(\delta(b))\varphi)(n(X)m(t))
 =\varphi(\delta(b)^{-1}n(X)m(t)\delta(b))
 =\varphi(n(b^{-1}X)m(t)).
\]
Substituting $X=bY$ in \eqref{eq:Bessel-period} proves the assertion.
\end{proof}

Put
\[
 C_{\Pi,\infty}
 =\frac{L_\infty(1/2,\Pi)
 L_\infty(1/2,\Pi\otimes\operatorname{sgn})}
 {L_\infty(1,\Pi,\Ad)L_\infty(1,\operatorname{sgn})}>0.
\]
Here $L_\infty(1,\operatorname{sgn})=\Gamma_\R(2)=\pi^{-1}$.
Since $\chi_{d,\infty}=\operatorname{sgn}$, this conversion is independent of $d$, and
\begin{equation}\label{eq:archimedean-normalization-check}
 C_{\Pi,\infty}\mathcal N_\infty
 =\frac{1}{\zeta_\R(2)\zeta_\R(4)},\qquad
 C_{\Pi,\infty}\mathcal N_\infty\zeta_\Q(2)\zeta_\Q(4)
 =\zeta(2)\zeta(4).
\end{equation}

For the local torus measures of
Subsection~\ref{subsec:scalar-Bessel-quotient}, the ideal class decomposition gives
\[
 \vol_{\mathrm{prod}}
 \bigl(\A^\times T_d(\Q)\backslash T_d(\A)\bigr)
 =\frac{2h(K_d)}{w(K_d)}.
\]
See also \cite[\S~8.2.3]{FurusawaMorimoto2024}.
Comparison with Tamagawa volume $2$ gives
\begin{equation}
 \label{eq:torus-measure-factor}
 dt_{\mathrm{Tam}}=\mathfrak c_d\prod_vdt_v,
 \qquad
 \mathfrak c_d=\frac{w(K_d)}{h(K_d)}.
\end{equation}
The finite class-number formula gives $\mathfrak c_d=2\pi/\big(\sqrt{|d|}L(1,\chi_d)\big)$. Substitution in \eqref{eq:L-weight} yields \begin{equation}\label{eq:L-weight-square-root}
 \mathscr L_d(\Pi)^{1/2}
 =C_\Pi |d|^{-1/4}
 \frac{|L(1/2,\Pi\otimes\chi_d)|^{1/2}}{L(1,\chi_d)},\qquad  C_\Pi=
 \left(\frac{2\pi\zeta_\Q(2)\zeta_\Q(4)
 |L(1/2,\Pi)|}{S_\Pi|L(1,\Pi,\Ad)|}\right)^{1/2}.
\end{equation}
Since the spin transfer $\Sym^3\pi_f$ is cuspidal, its isobaric
decomposition has one summand, so $\ell(\Pi)=1$ and $S_\Pi=2$;
see \cite[(1.18)--(1.19)]{FurusawaMorimoto2021}. We normalize the global period by
\begin{equation}
 \label{eq:normalized-period}
 \mathcal B^\natural_{d,b}(\varphi)=
 \frac{\mathcal B_{d,b}(\varphi)}
 {(\mathcal J_{\infty,d,b})^{1/2}}
 \quad(b>0),\qquad
 \mathcal B^\natural_{d,b}(\varphi)=0\quad(b\leq0).
\end{equation}
Here $\mathcal J_{\infty,d,b}>0$ for $b>0$
by Proposition~\ref{prop:real-factor}. We are now ready to prove Theorem~\ref{thm:L-bound}.

\begin{proof}[Proof of Theorem~\ref{thm:L-bound}]
Fix $d\equiv a_0\pmod{M_\varphi}$ and the factorizable vector
$\varphi$.  By \eqref{eq:normalized-period} and
Proposition~\ref{prop:real-factor}, it suffices to consider $b>0$ with
$\mathcal B_{d,b}(\varphi)\neq0$. In particular, by Lemma~\ref{lem:global-covariance},
\[
 \mathcal B_{d,1}\bigl(R(\delta(b))\varphi\bigr)
 =\mathcal B_{d,b}(\varphi)\ne0.
\]

Fix a nonzero finite tensor $\eta_{\mathrm{fin}}$ and an open compact subgroup
$K_{\mathrm{fin}}$ fixing $\eta_{\mathrm{fin}}$. For each smooth vector $u$ of $\Pi_\infty$,
put $F_u=u\otimes\eta_{\mathrm{fin}}$.
By \eqref{eq:compact_quotient} and Sobolev restriction at level $K_{\mathrm{fin}}$,
there are finitely many $D_i\in U(\mathfrak g_\C)$ with
\begin{equation}\label{eq:global-period-continuity}
 |\mathcal B_{d,1}(F_u)|
 \leq v_d\sup_{r\in[R_d]_Z}|F_u(r)|
 \ll_{d,K_{\mathrm{fin}}}\sum_i\|R(D_i)F_u\|
 \ll_{d,\eta_{\mathrm{fin}}}\sum_i\|d\Pi_\infty(D_i)u\|_\infty.
\end{equation}
The period is therefore smooth-continuous. Taking $\eta_{\mathrm{fin}}$ from
$R(\delta(b))\varphi$ and using $K$-finite density supplies a nonzero
$K$-finite period. Applying \cite[Theorem~1]{FurusawaMorimoto2021} gives, initially for $K$-finite
$\eta_\infty$ and $\eta=\eta_\infty\otimes\eta_{\mathrm{fin}}\ne0$,\begin{equation}
 \label{eq:refined-formula-one}
 \frac{|\mathcal B_{d,1}(\eta)|^2}
 {\langle\eta,\eta\rangle}
 =C_{\Pi,\infty}\frac{\mathfrak c_d}{S_\Pi}
 \frac{\zeta_\Q(2)\zeta_\Q(4)L(1/2,\Pi)
 L(1/2,\Pi\otimes\chi_d)}
 {L(1,\Pi,\Ad)L(1,\chi_d)}
 \prod_vJ_{v,d,1}(\eta_v).
\end{equation}

For general $\eta_\infty\ne0$, take $K$-finite $u_n\to\eta_\infty$
smoothly and set $\eta_n=u_n\otimes\eta_{\mathrm{fin}}$, with finite factors fixed.
By \eqref{eq:global-period-continuity} the periods converge, and
Lemma~\ref{lem:real-lowest-type} gives convergence of the real relative form.
Also,
\[
 \|\eta_n-\eta\|\ll_{\eta_{\mathrm{fin}}}\|u_n-\eta_\infty\|_\infty\longrightarrow0,
 \qquad \|\eta_n\|^2\longrightarrow\|\eta\|^2>0.
\]
Thus the normalized real factor converges, extending
\eqref{eq:refined-formula-one} to every nonzero decomposable smooth vector.

Apply it to
$R(\delta(b))\varphi=\otimes_v\Pi_v(\delta(b))\varphi_v$.
This is a smooth pure tensor, since $\delta(b)$ fixes $\varphi_p$ at
almost every prime. Its norm is unchanged and its period is
$\mathcal B_{d,b}(\varphi)$ by Lemma~\ref{lem:global-covariance}. Moreover, at every place,
$J_{v,d,1}\bigl(\Pi_v(\delta(b))\varphi_v\bigr) =|b|_v^3J_{v,d,b}(\varphi_v)$.
The product formula $\prod_v|b|_v^3=1$ cancels the Jacobians.
Dividing by
$\mathcal J_{\infty,d,b}>0$, and taking absolute values gives
\begin{equation*}
 \frac{|\mathcal B^\natural_{d,b}(\varphi)|^2}
 {\langle\varphi,\varphi\rangle}
 =C_{\Pi,\infty}\mathscr L_d(\Pi)
 \frac{J_{\infty,d,b}(\varphi_\infty)}
 {\mathcal J_{\infty,d,b}}
 \prod_{p<\infty}|J_{p,d,b}(\varphi_p)|.
\end{equation*}
For the original $\varphi_\infty\in W_\infty$, the real quotient lies
in $[0,1]$ by \eqref{eq:archimedean-comparison}.
The nonvanishing assumption supplies the denominator support condition
in Proposition~\ref{prop:finite-product}; its bound
$O_\varphi(H(b)^3)$ proves \eqref{eq:L-bound}.
\end{proof}
\subsection{Arithmetic estimates and CM coefficients}

Since $8\mid M_\varphi$, $a_0\equiv1\pmod4$ and
$\gcd(a_0,M_\varphi)=1$, the family corresponds to squarefree
$n=|d|\in[D,2D]$ with $n\equiv-a_0\pmod{M_\varphi}$.
M\"obius inversion gives
\begin{equation}\label{eq:family-count}
 \#\cD(D,a_0)
 =\sum_{\substack{r\leq\sqrt{2D}\\\gcd(r,M_\varphi)=1}}
 \mu(r)\left(\frac{D}{M_\varphi r^2}+O(1)\right)
 =\frac{D}{M_\varphi}
 \prod_{p\nmid M_\varphi}(1-p^{-2})+O_{M_\varphi}(D^{1/2}).
\end{equation}
In particular, $\#\mathcal D(D,a_0)\asymp_{M_\varphi}D$
for sufficiently large $D$.

\begin{proposition}\label{prop:arithmetic-moment}
Assume the GRH for
$L(s,\chi_d)$ and $L(s,\Sym^3\pi_f\otimes\chi_d)$
for every negative fundamental discriminant
$d\equiv a_0\pmod{M_\varphi}$. For every $\varepsilon>0$, we have
\begin{equation}\label{eq:arithmetic-pointwise}
 \mathscr L_d(\Pi)^{1/2}
 \ll_{\Pi,\varepsilon}|d|^{-1/4+\varepsilon}
 \qquad(d\equiv a_0\pmod{M_\varphi}).
\end{equation}

For the mean estimate, assume in addition GRH for
$L(s,\Sym^k\pi_f)$, $k=2,4,6$, and for both displayed twists at every negative fundamental discriminant $d\equiv a_0\pmod{M_0}$. Then, for every $\varepsilon>0$ and $D\geq3$, we have
\begin{equation}\label{eq:arithmetic-power-saving}
 \frac1D\sum_{d\in\mathcal D(D,a_0)}\mathscr L_d(\Pi)^{1/2}
 \ll_{\Pi,M_\varphi,\varepsilon}
 D^{-1/4}(\log D)^{-1/8+\varepsilon}.
\end{equation}
\end{proposition}

\begin{proof}
Since $\gcd(d,2M_\Pi)=1$, the spin twist has finite conductor
$M_\Pi|d|^4$ and fixed real parameters. Thus its analytic conductor
$\mathcal C_d$ at the central point satisfies
\begin{equation}\label{eq:twist-analytic-conductor}
 \mathcal C_d\asymp_\Pi |d|^4.
\end{equation}
Under the pointwise GRH hypotheses,
\cite[Theorem~2.1]{Chandee2009} and
\cite[Lemmas~2.3 and~2.5]{LamzouriLiSoundararajan2015} give
\[
 |L(1/2,\Pi\otimes\chi_d)|^{1/2}\ll_{\Pi,\tau}|d|^\tau,
 \qquad L(1,\chi_d)^{-1}\ll\log\log(3|d|)
 \quad(\tau>0).
\]
Here we take $x=(\log|d|)^2$ and use the Ramanujan bounds. Substitution in \eqref{eq:L-weight-square-root}, with
$\tau=\varepsilon/2$, proves \eqref{eq:arithmetic-pointwise};
bounded discriminants are absorbed in the constants.

For the mean, let $\mathcal D_0(D,a_0)$ have modulus $M_0$ in place
of $M_\varphi$. Then $\mathcal D(D,a_0)\subseteq\mathcal D_0(D,a_0)$,
and the required twist hypotheses hold throughout this larger family.
The lifts $\Sym^j\pi_f$, $j=2,4,6$, are cuspidal by
\cite[Theorem~A]{NewtonThorne2021}.
Let $N_f$ be the level of $f$ and let $a_\Pi(p^2)$ denote the
logarithmic Euler coefficient at every prime.
For $p\nmid N_f$, write $\lambda_j(p)=\lambda_{\Sym^j\pi_f}(p)$;
then
\[a_\Pi(p^2)=\sum_{i=1}^4\alpha_{i,p}^2\]
for the spin Satake
parameters $\alpha_{i,p}$, and
\begin{equation}\label{eq:arithmetic-satake-identities}
 \lambda_3(p)^2=1+\lambda_2(p)+\lambda_4(p)+\lambda_6(p),\qquad
 a_\Pi(p^2)=\lambda_6(p)-\lambda_4(p)+\lambda_2(p)-1.
\end{equation}
Fixed-lift GRH gives
$\sum_{p\leq t,\,p\nmid N_f}\lambda_j(p)/p=O_f(1)$ for $j=2,4,6$.
Since $\sum_{p\mid d}p^{-1}\ll1+\log\log\log(3D)$, Mertens' theorem
and partial summation yield
\[
 \begin{aligned}
 \sum_{\substack{y<p\leq z\\p\nmid dN_f}}\frac{\lambda_3(p)^2}{p}
 &=\log\frac{\log z}{\log y}
   +O_f\bigl(1+\log\log\log(3D)\bigr),\\
 \frac12\sum_{\substack{p\leq\sqrt x\\p\nmid d}}
 \frac{a_\Pi(p^2)w_x(p)}p
 &=-\frac12\log\log x
   +O_f\bigl(1+\log\log\log(3D)\bigr),
 \end{aligned}
\]
where $2\leq y\leq z$, $x\geq4$, and
$w_x(t)=t^{-2/\log x}(1-2\log t/\log x)$ decreases from at most $1$
to $0$ on $[2,\sqrt x]$.
Ramified primes and higher prime powers contribute $O_f(1)$.
These give the single-form variance and square-term inputs in
\cite[Lemma~4.3 and (5.1)--(5.2)]{HuaMiao2025}, without RH for $\zeta$.
Thus \cite[Theorem~3.2]{HuaMiao2025} applies with $m=1$, $\ell_1=1/2$,
whose exponent is $\ell_1(\ell_1-1)/2=-1/8$.

In the large-$V$ case of \cite[Lemma~4.5]{HuaMiao2025}, take
$r=\lfloor\epsilon V/10\rfloor$, where $\epsilon>0$ is fixed
and sufficiently small. Markov's inequality and Stirling's formula
bound the number of discriminants for which either prime sum
exceeds its threshold by
\[
 O_{\Pi,\epsilon}\!\left(D e^{-c_\epsilon V\log V}\right),
 \qquad c_\epsilon>0,
\]
uniformly for sufficiently large $D$. This tail bound suffices
for the half-moment estimate.

Self-duality and GRH imply $L(1/2,\Pi\otimes\chi_d)\geq0$. Hence, for every $\tau>0$,
\begin{align*}
 \frac1D\sum_{d\in\mathcal D(D,a_0)}\mathscr L_d(\Pi)^{1/2}
 &\ll_{\Pi,M_\varphi}
 D^{-1/4}\log\log(3D)\,
 \frac1D\sum_{d\in\mathcal D_0(D,a_0)}
 L(1/2,\Pi\otimes\chi_d)^{1/2}\\
 &\ll_{\Pi,M_\varphi,\tau}
 D^{-1/4}\log\log(3D)(\log D)^{-1/8+\tau}.
\end{align*}
Taking $\tau=\varepsilon/2$ and absorbing $\log\log(3D)$ proves
\eqref{eq:arithmetic-power-saving}.
By \eqref{eq:family-count}, the same bound holds with
$1/\#\mathcal D(D,a_0)$ for sufficiently large $D$.
\end{proof}

\begin{proposition}
\label{prop:uniform-character}
Under the pointwise hypotheses of
Proposition~\ref{prop:arithmetic-moment}, for every $b\in\Q_{>0}$, $\varepsilon>0$, and negative fundamental discriminant $d\equiv a_0\pmod{M_\varphi}$, we have
\begin{align*}
 |\mathcal C_{d,b}(\varphi)|
 \ll_{\varphi,\varepsilon}H(b)^{3/2}
 |d|^{-1/4+\varepsilon},
 \qquad
 |\mathcal B^\natural_{d,b}(\varphi)|
 \ll_{\varphi,\varepsilon}H(b)^{3/2}
 |d|^{-1/4+\varepsilon}.
\end{align*}
Under the additional mean hypotheses, for every $D\geq3$, we have
\begin{align}
 \frac1D\sum_{d\in\mathcal D(D,a_0)}|\mathcal C_{d,b}(\varphi)|
 &\ll_{\varphi,\varepsilon}
 H(b)^{3/2} D^{-1/4}(\log D)^{-1/8+\varepsilon},
 \label{eq:CM-average-mean}\\
 \frac1D\sum_{d\in\mathcal D(D,a_0)}
 |\mathcal B^\natural_{d,b}(\varphi)|
 &\ll_{\varphi,\varepsilon}
 H(b)^{3/2} D^{-1/4}(\log D)^{-1/8+\varepsilon}.
 \label{eq:uniform-character}
\end{align}
All coefficients vanish for $b\leq0$.
\end{proposition}

\begin{proof}
For $b>0$, \eqref{eq:CM-average-Fourier} and \eqref{eq:normalized-period} give
\[
 \mathcal C_{d,b}(\varphi)
 =\frac12\bigl(\mathcal J_{\infty,d,b}\bigr)^{1/2}
   \mathcal B^\natural_{d,b}(\varphi).
\]
Theorem~\ref{thm:L-bound} and
Corollary~\ref{cor:real-uniform-bound} therefore give
the unconditional bound
\begin{equation}\label{eq:coefficient-arithmetic-input}
 |\mathcal C_{d,b}(\varphi)|+|\mathcal B^\natural_{d,b}(\varphi)|
 \ll_\varphi H(b)^{3/2}\mathscr L_d(\Pi)^{1/2}.
\end{equation}
Substituting the pointwise or mean estimate of
Proposition~\ref{prop:arithmetic-moment} proves the bounds for $b>0$.
The vanishing for $b\leq0$ follows from
Proposition~\ref{prop:zero-and-negative} and
\eqref{eq:normalized-period}.
\end{proof}

\begin{theorem}
\label{thm:raw-exponential}
With $f,\Pi,\varphi$ as in Subsection~\ref{subsec:standing-data},
put $Q=Q_\varphi$.
For every $\varepsilon>0$ and every negative fundamental discriminant
$d\equiv a_0\pmod{M_\varphi}$, we have
\begin{equation*}
 \|\mathscr A_d(\varphi)\|_{\mathrm{TV}}
 \ll_{\varphi,\varepsilon}
 |d|^{\kappa-1/2+\varepsilon}
 \exp\!\left(-\frac{\pi|d|}{2Q}\right).
\end{equation*}
Consequently, for every bounded Borel function $\Phi$ on
$\mathscr X=\Q\backslash\A$,
\[
 |\langle\mathscr A_d(\varphi),\Phi\rangle|
 \ll_{\varphi,\varepsilon}
 \|\Phi\|_\infty |d|^{\kappa-1/2+\varepsilon}
 \exp\!\left(-\frac{\pi|d|}{2Q}\right).
\]
\end{theorem}

\begin{proof}
Put $D=|d|$. By Lemma~\ref{lem:denominator} and
Proposition~\ref{prop:zero-and-negative},
$\mathcal C_{d,b}(\varphi)=0$ unless $b=j/Q$ with $j\geq1$.
If $\sum_{j\geq1}|\mathcal C_{d,j/Q}(\varphi)|<\infty$, absolute convergence
and Fourier uniqueness give
\[
 \mathscr A_d(\varphi)=a_d(x)\,dx,\qquad
 a_d(x)=\sum_{j\geq1}
 \mathcal C_{d,j/Q}(\varphi)\psi(jx/Q).
\]
Hence
\[
 \|\mathscr A_d(\varphi)\|_{\mathrm{TV}}
 =\int_{\mathscr X}|a_d(x)|\,dx
 \leq\sum_{j\geq1}|\mathcal C_{d,j/Q}(\varphi)|.
\]
It therefore suffices to bound
$\sum_{j\geq1}|\mathcal C_{d,j/Q}(\varphi)|$.

For $b>0$, the period normalization and
Theorem~\ref{thm:L-bound} give
\[
 |\mathcal C_{d,b}(\varphi)|
 =\frac12
   \bigl(\mathcal J_{\infty,d,b}\bigr)^{1/2}
   |\mathcal B^\natural_{d,b}(\varphi)|
 \ll_\varphi
   \bigl(\mathcal J_{\infty,d,b}\bigr)^{1/2}
   H(b)^{3/2}\mathscr L_d(\Pi)^{1/2}.
\]
The conductor estimate \eqref{eq:twist-analytic-conductor},
the convexity bound, and the class-number formula give
\[
 |L(1/2,\Pi\otimes\chi_d)|^{1/2}
 \ll_{\Pi,\varepsilon}D^{1/2+\varepsilon},
 \qquad
 L(1,\chi_d)^{-1}
 =\frac{w(K_d)\sqrt D}{2\pi h(K_d)}
 \ll D^{1/2},
\]
where $h(K_d)\geq1$ and $w(K_d)\leq6$.
Thus \eqref{eq:L-weight-square-root} yields
\[
 \mathscr L_d(\Pi)^{1/2}
 \ll_{\Pi,\varepsilon}D^{3/4+\varepsilon}.
\]
Combining this with \eqref{eq:real-absolute-bound}, we obtain
\[
 |\mathcal C_{d,b}(\varphi)|
 \ll_{\varphi,\varepsilon}
 b^{(3\kappa-3)/2}H(b)^{3/2}
 D^{\kappa-1/2+\varepsilon}e^{-\pi bD/2},
\]
since $\tr(S_d)\geq D/4$.
For $b=j/Q$ with $j\geq1$, we have
$H(b)\leq\max\{j,Q\}\leq Qj$. Therefore,
\[
 \begin{aligned}
 \sum_{j\geq1}|\mathcal C_{d,j/Q}(\varphi)|
 &\ll_{\varphi,\varepsilon}
 D^{\kappa-1/2+\varepsilon}e^{-\pi D/(2Q)}
 \sum_{j\geq1}j^{3\kappa/2}
 e^{-\pi(j-1)D/(2Q)}\\
 &\ll_{\varphi,\varepsilon}
 D^{\kappa-1/2+\varepsilon}e^{-\pi D/(2Q)},
 \end{aligned}
\]
where the remaining series is bounded uniformly for $D\geq3$. This completes the proof.
\end{proof}

\section{Fourier summation}
\label{sec:Fourier}

Recall $\mathscr X=\Q\backslash\A$.
For an open compact additive subgroup $U\subset\A_{\mathrm{fin}}$, let $C^\infty(\mathscr X)^U$ denote the $U$-invariant functions smooth in
the real direction, with Fr\'echet seminorms
\[
 p_J(\Phi)=\sum_{j=0}^J
 \|\partial_\infty^j\Phi\|_{L^\infty(\mathscr X)},\qquad J\geq0.
\]
Give $C^\infty_{\mathrm{fl}}(\mathscr X)=\bigcup_U C^\infty(\mathscr X)^U$
the locally convex inductive-limit topology, with cofinal levels
$n!\widehat\Z$. A distribution is a complex-linear functional continuous
on each finite-level Fr\'echet space. Pairings are complex-linear in the test function.
For $\Phi\in C^\infty_{\mathrm{fl}}(\mathscr X)$, let
\[
 \widehat\Phi(b)=\int_{\mathscr X}\Phi(x)\psi(-bx)\,dx,
 \qquad b\in\Q.
\]
The characters $x\mapsto\psi(bx)$ identify the Pontryagin dual of $\mathscr X$ with $\Q$, and the finite part $\psi_{\mathrm{fin}}=\prod_{p<\infty}\psi_p$ has conductor $\widehat\Z$.
\begin{lemma}
\label{lem:Fourier-decay}
For every $\Phi\in C^\infty_{\mathrm{fl}}(\mathscr X)$, there is an integer $Q_\Phi\geq1$ such that $\widehat\Phi(m/n)=0$ unless $n\mid Q_\Phi$, where $m/n$ is in lowest terms and $n>0$.
For every $B\geq0$, we have
\begin{equation}
 \label{eq:Fourier-decay}
 |\widehat\Phi(b)|\ll_{\Phi,B}(1+H(b))^{-B}.
\end{equation}
For every $A,L\geq0$ and $T\geq1$,
\begin{equation}
 \label{eq:Fourier-tail}
 \sum_{H(b)>T}H(b)^A|\widehat\Phi(b)|
 \ll_{\Phi,A,L}T^{-L}.
\end{equation}
More generally, for every integer $Q\geq1$, every $A\geq0$, and
every integer $J>A+1$,
\begin{equation}\label{eq:fixed-lattice-summation}
 \sum_{b\in Q^{-1}\Z}H(b)^A|\widehat\Phi(-b)|
 \ll_{Q,A,J}p_J(\Phi),
 \qquad \Phi\in C^\infty_{\mathrm{fl}}(\mathscr X),
\end{equation}
independently of the finite level of $\Phi$. In particular, at a fixed
finite level $U$, the sum over all $b\in\Q$ is
$O_{U,A,J}(p_J(\Phi))$.
\end{lemma}

\begin{proof}
Choose $U$ fixing $\Phi$ and $Q_U\widehat\Z\subset U$.
Translation gives $\widehat\Phi(b)=\psi_{\mathrm{fin}}(bu)\widehat\Phi(b)$ for $u\in U$.
The conductor of $\psi_{\mathrm{fin}}$ forces $bQ_U\in\Z$ whenever
$\widehat\Phi(b)\ne0$; take $Q_\Phi=Q_U$.

Integration by parts gives $(2\pi ib)^j\widehat\Phi(b)
=\widehat{\partial_\infty^j\Phi}(b)$ for $b\ne0$, hence
\[
|\widehat\Phi(b)|\ll_J(1+|b|_\infty)^{-J}p_J(\Phi).
\]
For $b\in Q^{-1}\Z$, cancellation in $b=m/Q$ gives
$H(b)\leq Q\max\{1,|b|_\infty\}$. Hence
\[
 \sum_{b\in Q^{-1}\Z}H(b)^A|\widehat\Phi(-b)|
 \ll_{Q,A,J}p_J(\Phi)
 \sum_{m\in\Z}(1+|m|)^{A-J}.
\]
The series converges for $J>A+1$, proving
\eqref{eq:fixed-lattice-summation}. Taking $Q=Q_U$ and using the
Fourier support proves the fixed-level assertion and
\eqref{eq:Fourier-decay}.
For $J>A+L+1$, the same weighted estimate gives directly
\[
 \sum_{H(b)>T}H(b)^A|\widehat\Phi(b)|
 \leq T^{-L}\sum_{b\in\Q}H(b)^{A+L}|\widehat\Phi(b)|
 \ll_{U,A,L,J}T^{-L}p_J(\Phi),
\]
which proves \eqref{eq:Fourier-tail}.
\end{proof}

Lemma~\ref{lem:Fourier-decay} and Fourier uniqueness give the
absolutely and uniformly convergent Fourier expansion $\Phi(x)=\sum_b\widehat\Phi(b)\psi(bx)$. Thus, for a finite complex measure $\mu$,
\[
 \langle\mu,\Phi\rangle
 =\sum_{b\in\Q}\widehat\mu(b)\widehat\Phi(-b),
 \qquad
 \widehat\mu(b)=\int_{\mathscr X}\psi(-bx)\,d\mu(x).
\]
In particular, $\mathscr A_d(\varphi)$ has Fourier coefficients
$\mathcal C_{d,b}(\varphi)=\mathcal B_{d,b}(\varphi)/2$.

\begin{definition}
\label{def:normalized-distribution}
For every negative fundamental discriminant
$d\equiv a_0\pmod{M_\varphi}$, define
\begin{equation}
 \label{eq:normalized-distribution}
 \left\langle\mathscr P_d^\natural(\varphi),\Phi\right\rangle
 =\sum_{b\in\Q}
 \mathcal B^\natural_{d,b}(\varphi)\widehat\Phi(-b),
 \qquad\Phi\in C^\infty_{\mathrm{fl}}(\mathscr X).
\end{equation}
Here the normalized global period $\mathcal B^\natural_{d,b}(\varphi)$ is defined as in \eqref{eq:normalized-period}.
\end{definition}

\begin{proposition}
\label{prop:finite-order}
The series \eqref{eq:normalized-distribution} converges absolutely and defines a distribution.  For every $\Phi\in C^\infty_{\mathrm{fl}}(\mathscr X)$, we have
\begin{equation}
 \label{eq:finite-order-seminorm}
 \left|\left\langle
 \mathscr P_d^\natural(\varphi),\Phi
 \right\rangle\right|
 \ll_\varphi
 \mathscr L_d(\Pi)^{1/2}p_3(\Phi).
\end{equation}
The same estimate holds for $\mathscr A_d(\varphi)$, with implied constants independent of the finite level of $\Phi$.
\end{proposition}

\begin{proof}
By \eqref{eq:coefficient-arithmetic-input}, both coefficient families
are bounded by $O_\varphi(H(b)^{3/2}\mathscr L_d(\Pi)^{1/2})$.
Lemma~\ref{lem:denominator} and
Proposition~\ref{prop:zero-and-negative} restrict their support to
$\Q_{>0}\cap Q_\varphi^{-1}\Z$.
Apply \eqref{eq:fixed-lattice-summation} with $Q=Q_\varphi$,
$A=3/2$, and $J=3$. This gives
\[
 \sum_{\substack{b>0\\b\in Q_\varphi^{-1}\Z}}
 H(b)^{3/2}|\widehat\Phi(-b)|\ll_\varphi p_3(\Phi),
\]
independently of the finite level of $\Phi$.
Together with \eqref{eq:coefficient-arithmetic-input}, this proves
absolute convergence and \eqref{eq:finite-order-seminorm} for both
families.  Continuity on each finite-level space, and hence on the
stated locally convex inductive limit, follows from the same estimate.
\end{proof}

\begin{theorem}\label{thm:distribution}
With $f,\Pi,\varphi$ as in Subsection~\ref{subsec:standing-data},
assume the GRH for
$L(s,\chi_d)$ and $L(s,\Sym^3\pi_f\otimes\chi_d)$
for every negative fundamental discriminant
$d\equiv a_0\pmod{M_\varphi}$.
For every $\varepsilon>0$, every such $d$, and every $\Phi\in C^\infty_{\mathrm{fl}}(\mathscr X)$, we have
\begin{equation}\label{eq:pointwise-distribution}
 \left|\left\langle\mathscr P_d^\natural(\varphi),\Phi\right\rangle\right|
 \ll_{\varphi,\varepsilon}
 p_3(\Phi)|d|^{-1/4+\varepsilon}.
\end{equation}
For the mean estimate, assume in addition GRH for
$L(s,\Sym^k\pi_f)$, $k=2,4,6$, and for both displayed twists at every negative fundamental discriminant $d\equiv a_0\pmod{M_0}$.
Then, for every $\varepsilon>0$ and $D\geq3$, we have
\begin{equation}\label{eq:distribution}
 \frac1D\sum_{d\in\mathcal D(D,a_0)}
 \left|\left\langle\mathscr P_d^\natural(\varphi),\Phi\right\rangle\right|
 \ll_{\varphi,\varepsilon}
 p_3(\Phi)D^{-1/4}(\log D)^{-1/8+\varepsilon}.
\end{equation}
Both estimates are uniform over all finite levels of $\Phi$.
\end{theorem}
\begin{proof}
Substituting \eqref{eq:arithmetic-pointwise} into
\eqref{eq:finite-order-seminorm} proves \eqref{eq:pointwise-distribution}.
Averaging \eqref{eq:finite-order-seminorm} over $\mathcal D(D,a_0)$ and
applying \eqref{eq:arithmetic-power-saving} proves \eqref{eq:distribution}.
All constants are uniform on $p_3$-bounded subsets of $C^\infty_{\mathrm{fl}}(\mathscr X)$, independently of the finite level.
\end{proof}
\begin{proof}[Proof of Theorem~\ref{thm:pointwise-distribution}]
For fixed $\Phi$, \eqref{eq:pointwise-distribution} gives the stated
bound, and taking $0<\varepsilon<1/4$ yields the limit.
\end{proof}
As an application, we approximate the distribution pairings by finite
Fourier sums, with a quantitative bound for the mean truncation error.
\begin{corollary}
\label{cor:truncation}
Assume the mean hypotheses of
Proposition~\ref{prop:arithmetic-moment}, and fix
$\Phi\in C^\infty_{\mathrm{fl}}(\mathscr X)$.
For every $\varepsilon>0$, $L\geq0$, $D\geq3$, and $T\geq1$, we have
\begin{equation*}
 \frac1D\sum_{d\in\mathcal D(D,a_0)}
 \left|
 \left\langle\mathscr A_d(\varphi),\Phi\right\rangle
 -\sum_{\substack{b\in\Q\\H(b)\leq T}}
 \mathcal C_{d,b}(\varphi)\widehat\Phi(-b)
 \right|
 \ll_{\varphi,\Phi,\varepsilon,L}
 T^{-L}D^{-1/4}(\log D)^{-1/8+\varepsilon}.
\end{equation*}
The same bound holds with $\mathscr A_d$ and $\mathcal C_{d,b}$ replaced by $\mathscr P_d^\natural$ and $\mathcal B^\natural_{d,b}$.
Every nonzero term in either truncated sum has $b=m/n$ in lowest terms with
\[
 n\mid\gcd(Q_\Phi,Q_\varphi),
 \qquad 1\leq m\leq T,\qquad 1\leq n\leq T.
\]
\end{corollary}

\begin{proof}
By absolute convergence, the difference is the Fourier tail over
$H(b)>T$.  The triangle inequality, Tonelli's theorem for the resulting
nonnegative series, and \eqref{eq:CM-average-mean} give
\begin{equation*}
 \sum_{H(b)>T}|\widehat\Phi(-b)|
 \frac1D\sum_{d\in\mathcal D(D,a_0)}|\mathcal C_{d,b}(\varphi)|
 \ll_{\varphi,\varepsilon}
 D^{-1/4}(\log D)^{-1/8+\varepsilon}
 \sum_{H(b)>T}H(b)^{3/2}|\widehat\Phi(-b)|.
\end{equation*}
Apply \eqref{eq:Fourier-tail} with $A=3/2$.
For the normalized distribution use \eqref{eq:uniform-character} instead.
The support conditions give $b=m/n>0$ and
$n\mid\gcd(Q_\Phi,Q_\varphi)$; $H(b)\le T$ gives $m,n\le T$.
\end{proof}

\section*{Acknowledgements}
S.H. would like to thank Chung Pang Mok for the constant encouragement.
S.Y. would like to thank Ralf Schmidt for his helpful discussions.
S.H. is supported by the National Natural Science Foundation of China (No. 12601008) and the China Postdoctoral Science Foundation (No. 2026M793348).
S.Y. is supported by the National Natural Science Foundation of China (Nos. 12301016, 12471187) and the Fundamental Research Funds for the Central Universities (No. 20720260221).

\section*{AI Use Disclosure}
The research direction and initial mathematical framework were developed by the authors.
The pointwise distribution estimate in Theorem~\ref{thm:pointwise-distribution}
is deduced from the polynomial Bessel-period bound in
Theorem~\ref{thm:L-bound}, together with arithmetic estimates under GRH
and Fourier summation. Before using generative AI, the authors had obtained by explicit computation the $j=1$ case of Lemma~\ref{lem:good-place} and the real Bessel factor in the scalar-valued case (Subsection~\ref{subsec:real-factor}). These calculations supplied key local ingredients for the polynomial bound in the scalar-valued setting with $b$ a positive squarefree integer. Generative AI was subsequently used interactively in the mathematical exploration. The authors independently verified the final arguments and take full responsibility for the paper.

\end{document}